\documentclass[12pt, a4paper, oneside]{amsart}
\usepackage[a4paper, margin = 2cm]{geometry}
\usepackage[english]{babel}
\usepackage[T2A]{fontenc}
\usepackage[utf8]{inputenc} 
\usepackage{amsfonts}
\usepackage{amssymb, amsthm, amscd}
\usepackage{amsmath}
\usepackage{mathtools}
\usepackage{needspace}
\usepackage{etoolbox}
\usepackage{lipsum}
\usepackage{comment}
\usepackage{cmap}
\usepackage[pdftex]{graphicx}
\usepackage[unicode]{hyperref}
\usepackage[matrix,arrow,curve]{xy}
\usepackage[usenames,dvipsnames]{xcolor}
\usepackage{colortbl}
\usepackage{textcomp}
\usepackage{cite}
\usepackage{euscript}
\usepackage{tikz}

\DeclareFontFamily{U}{mathx}{}
\DeclareFontShape{U}{mathx}{m}{n}{<-> mathx10}{}
\DeclareSymbolFont{mathx}{U}{mathx}{m}{n}
\DeclareMathAccent{\widehat}{0}{mathx}{"70}
\DeclareMathAccent{\widecheck}{0}{mathx}{"71}

\makeatletter
\def\@settitle{\begin{center}
		\baselineskip14\p@\relax
		\bfseries
		\LARGE
		\@title
	\end{center}
}
\makeatother

\newtheorem{theorem}{Theorem}[section]

\newtheorem{lemma}[theorem]{Lemma}
\newtheorem{corollary}[theorem]{Corollary}

\newtheorem{mainth}{Theorem}

\renewcommand{\geq}{\geqslant}
\renewcommand{\leq}{\leqslant}
\newcommand{\R}{\mathbb{R}}
\newcommand{\N}{\mathbb{N}}
\newcommand{\Z}{\mathbb{Z}}

\newcommand{\eps}{\varepsilon}

\newcommand{\const}{\mathrm{const}}
\newcommand{\CC}{\mathbb{C}}

\DeclareMathOperator{\Span}{Span}
\DeclareMathOperator{\Gram}{Gram}

\DeclareMathOperator{\dist}{dist}
\DeclareMathOperator{\sinc}{sinc}

\newenvironment{enumerate*}%
  {\begin{enumerate}%
    \setlength{\itemsep}{1pt}%
    \setlength{\parskip}{1pt}}%
  {\end{enumerate}}

\makeatletter
\newcommand*{\@old@slash}{}\let\@old@slash\slash
\def\slash{\relax\ifmmode\delimiter"502F30E\mathopen{}\else\@old@slash\fi}
\makeatother

\title{Schauder exponential basis which is not Riesz}
\author{Andrei V. Semenov}

\thanks{The work is supported by RSF grant 26-11-00157}
\keywords{Functional analysis, complex analysis, time-frequency analysis, basis, exponentials}

\address{Andrei V. Semenov
Department of Mathematics and Computer Science, Saint Petersburg State University\\
Russia, 199178, Saint Petersburg, 14th line of Vasilievsky island, 29}

\email{asemenov.spb.56@gmail.com}

\begin{document}

\begin{abstract}
We show that there exists a Schauder basis $\{e^{i \lambda_n t}\}_{n \in \N}$ with real frequencies in $L^2(-\pi,\pi)$ which is not a Riesz basis.
\end{abstract}

\maketitle

\section{Introduction}
Let $H$ be a separable infinite-dimensional Hilbert space. A system $\{e_n\}_{n \in \Z} \subset H$ is a Riesz basis if it is the image of an orthonormal basis under a bounded invertible operator on $H$. Equivalently, its linear span is dense in $H$ and there exist $A, B>0$ such that for any finitely supported sequence $\{c_n\} \in \ell^2(\Z)$ one has
$$A\sum_n|c_n|^2\leq \left\| \sum_{n \in \Z} c_n e_n\right\|_H^2 \leq B\sum_n |c_n|^2.$$
In this case for every $x \in H$ there is a unique expansion $x = \sum_n a_n(x) e_n$ converging unconditionally. \par

A system $\{e_n\}_{n \in \Z} \subset H$ is a Schauder basis if for a fixed enumeration of $\Z$ every $x \in H$ has a unique expansion $x = \sum_{n \in \Z} a_n(x) e_n$ converging in norm in that order. It is clear that every Riesz basis is a Schauder basis in every enumeration. The converse is false: a classical example due to Babenko is the system $\{e^{int}\}_{n\in\Z}$ in weighted $L^2((-\pi,\pi); |t|^\alpha)$ for $0<|\alpha|<1/2$ (see \cite{Babenko1948}). More generally, the Schauder basis property of the trigonometric system in weighted $L^2((-\pi,\pi); w)$ can be characterized by Muckenhoupt's $A_2$-condition (see \cite{HuntMuckenhouptWheeden1973}).\par

For exponential systems in the unweighted space $L^2(-\pi,\pi)$, the long-standing question is whether every Schauder basis $\{e^{i\lambda_nt}\}_{n\in\Z}$ with real frequencies must also be a Riesz basis. This question is discussed in Young's monograph (see ~\cite[p.~165]{Young2001}). Heil and Kutyniok in ~\cite[p.~573]{HeilKutyniok2008} explicitly noted that no such example was known and that its existence had not been ruled out. The question is also stated by Mitkovski in~\cite[Section~6.2.1, Problem~2]{Mitkovski2010}. Marchenko recorded the existence problem as open in 2024 (see~\cite[p.~753]{Marchenko2024}). Related questions concerning basis properties and stability were studied by Nakamura (see~\cite{Nakamura2007,Nakamura2009}). Connections between Riesz bases, Muckenhoupt conditions and complete interpolation in Paley--Wiener spaces are studied by Lyubarskii and Seip (see ~\cite{LyubarskiiSeip1997}).\par

The distinction between different basis properties of exponential systems arises naturally. For example, one may ask whether $L^2(D)$ admits a basis of exponentials with real frequencies for $D\subset \R^d$ of finite measure.  For convex bodies, Lev and Matolcsi proved that an orthogonal exponential basis exists if and only if the body tiles space by translations (see~\cite{LevMatolcsi2022}). Allowing Riesz bases enlarges the class of admissible domains: Kozma and Nitzan showed that every finite union of bounded intervals admits a Riesz basis of exponentials (see~\cite{KozmaNitzan2015}). \par 

On the other hand, a recent preprint of Ortega-Cerd\`a (see~\cite{OrtegaCerda2026}) establishes that balls in $\R^d$ for $d\geq2$ and nondegenerate planar triangles admit no Riesz basis of exponentials. These nonexistence results do not by themselves rule out Schauder bases.\par

In this article we resolve the question for $L^2(-\pi,\pi)$: there exists a Schauder basis of exponentials with real frequencies in $L^2(-\pi,\pi)$ which is not a Riesz basis. Moreover, we develop an algorithm for constructing such bases.

\section{Main result}
Consider $H = L^2(-\pi,\pi)$ with the usual norm and inner product
$$\|f\|_2^2 := \int_{-\pi}^{\pi} |f(t)|^2 \, dt  \quad \text{ and } \quad \langle f,g \rangle := \int_{-\pi}^{\pi} f(t)  \overline{g(t)} \, dt.$$

\begin{mainth}\label{th:main}
There exists a sequence $\{\lambda_n\}_{n \in \N} \subset \R$ such that the system $\{e^{i \lambda_n t}\}_{n \in \N}$ is a Schauder basis in $L^2(-\pi,\pi)$ but not a Riesz basis.
\end{mainth}

The main idea behind the construction can be informally stated as follows:
\begin{enumerate}
\item Take the standard basis $\{e^{int}\}_{n \in \Z} \subset L^2(-\pi,\pi)$ and consider a partition $\N = \bigcup_{k \geq 1} I_k$ into disjoint finite intervals $I_k$. Choose some $J_k \subset I_k$ for each $k$.
\item At the $k$-th stage remove every $e^{int}$ with $n \in J_k$ and insert $e^{i \lambda_{k,n} t}$ in its place, where $\lambda_{k,n} \in \Lambda_k \subset \R$. Within each $\Lambda_k$ use the order induced from $J_k$.
\item On each stage keep the norms of the coordinate projections onto initial segments of the insertions uniformly bounded. Make the smallest eigenvalues of the Gram matrices of the inserted exponentials tend to zero.
\item List $I_k\setminus J_k$ in increasing order. Let $\Lambda = \{\lambda_n \}_{n=1}^{\infty}$ be the final set of frequencies and enumerate it as follows:
$$0, \ \Lambda_1, \  I_1 \setminus J_1, \ -1, \ \Lambda_2, \  I_2 \setminus J_2,\ -2, \ \Lambda_3, \ I_3 \setminus J_3, \ -3, \ldots$$
\item Form a block system $\{y_{j}\}_{j \geq 1}$ by replacing $\{e^{i \lambda_{k,n} t}\}_{n \in J_k}$ with their orthogonal projections onto frequencies inside $I_k$. This system is a Schauder basis by orthogonality and the uniform projection bounds. 
\item Take each $I_k$ sufficiently large to make the tails of the series $e^{i\lambda_{k,n} t}= \sum_{l \in \Z} c_l(\lambda_{k,n}) e^{ilt}$ small. This leads to a bounded invertible operator $T \colon y_j \mapsto e^{i \lambda_{j} t}$. Thus the final system is a Schauder basis.
\item Finally, the system $\{e^{i \lambda_{k} t}\}_{\lambda_k \in \Lambda}$ has zero lower Riesz bound by the Gram matrix estimates.
\end{enumerate}

Section~\ref{sect:Alg} gives the precise algorithm and proves its correctness. A concrete example of such a basis is provided in Section~\ref{sect:Example}. The required estimates are proved in Section~\ref{sect:Tech}.

\subsection{Notation}
We define $\sinc x=\frac{\sin\pi x}{\pi x}$ and put $\sinc 0 = 1$. For an invertible matrix $M$ write $M^{-*}=(M^*)^{-1}$. For two Hermitian matrices $A,B$  of the same size we write $A \preceq B$ if for any vector $c$ we have $c^*Ac\leq c^*Bc$. By $|E|$ we denote the Lebesgue measure of $E\subset \R$. By $\#E$ we denote the number of elements of a finite set $E$.

\section{The Algorithm} \label{sect:Alg}
In this section we describe a construction based on successive finite replacements in the standard basis. \par
Put $e_z(t)=e^{izt}$ and $\widehat f(k)=\frac{\langle f,e_k\rangle}{2\pi}$. For $E\subset \Z$ we define $Q_E$ to be the orthogonal projection onto $H_E=\overline{\Span}\{e_k \mid k\in E\}$. Since $\{e_k\}_{k \in \Z}$ forms the standard basis, $Q_E$ is the corresponding coordinate projection and we have 
\begin{equation}\label{def:Q}
Q_E f=\sum_{k\in E} \widehat f(k) e_k.
\end{equation}

\subsection{One replacement}\label{sect:finrep} In this subsection we consider only one finite replacement of the standard basis. Fix $m \in \N$ and let $j_1<\cdots< j_m$ be integers. Let $\xi_1, \ldots, \xi_m \in \R$ be distinct. Put $J = \{j_1,\ldots, j_m\}$ and define $S \colon c \mapsto \sum_{l=1}^m c_l e_{\xi_l}$, $F \colon f \mapsto \frac1{2\pi}\bigl(\langle f, e_{j_l}\rangle\bigr)_{l=1}^m$, and $M=FS$. 

\begin{lemma}\label{lem:finrep}
The following conditions are equivalent:
\begin{enumerate}
\item $\det M\ne 0$,
\item the system $\{e_{\xi_l}\}_{l=1}^m \cup \{e_l\}_{l\in\Z\setminus J}$ is a Riesz basis in $L^2(-\pi,\pi)$,
\item one has $\{\xi_1,\ldots,\xi_m\}\cap \Z \subset J$.
\end{enumerate}
\end{lemma}
\begin{proof}
We first prove the equivalence of (1) and (2). Define the operators $T:\CC^m\oplus\ell^2(\Z\setminus J)\longrightarrow L^2(-\pi,\pi)$ and $U:L^2(-\pi,\pi) \longrightarrow\mathbb C^m\oplus\ell^2(\mathbb Z\setminus J)$ by the rule 
$$T(c,d)=\sum_{j=1}^m c_j e_{\xi_j}+\sum_{l\notin J}d_l e_l \quad \text{and} \quad Uf=\left(\{\widehat f(j_l)\}_{l=1}^m, \{\widehat f(l)\}_{l\in\mathbb Z\setminus J}\right).$$
Note that the operator $TU \colon L^2(-\pi,\pi) \to L^2(-\pi, \pi)$ maps the standard basis $\{e^{ilt}\}_{l \in \Z}$ to the system $\{e^{i \xi_l t}\}_{l=1}^m \cup \{e^{ilt}\}_{l \not\in J}$. We need to show that $TU$ has a bounded inverse. \par
It is clear that $U^{-1} \colon (c,d) \mapsto \sum_{l=1}^m c_l e^{i j_l t} + \sum_{l \not\in J} d_l e^{ilt}$ is bounded, since $j_l\in \Z$ and
$$\|U^{-1}(c,d)\|^2 = 2\pi \left(\sum_l |c_l|^2 +\sum_{l \not\in J} |d_l|^2\right) = 2\pi \|(c,d)\|^2.$$
Hence, $U$ has a bounded inverse and $\|Uf\|_2^2=\frac1{2\pi}\|f\|_2^2$. Now observe that
$$UT=\begin{pmatrix}M&0\\V&1\end{pmatrix},\quad \text{where} \quad V \colon c\mapsto \left\{ (Vc)_l=\sum_{j=1}^m c_j\sinc(\xi_j-l) \right\}_{l \in \Z \setminus J}.$$
The operator $V$ is bounded. Indeed, we have $\widehat{Sc}(l) = \sum_{k=1}^m  c_k \sinc(\xi_k-l)$ and
$$\|Vc\|^2_2 = \sum_{l \not\in J} \left| \sum_{k=1}^m c_k \sinc(\xi_k-l)\right|^2 = \sum_{l \not\in J} |\widehat{Sc}(l)|^2 \leq \sum_{l\in \Z}  |\widehat{Sc}(l)|^2  = {{1} \over {2\pi}} \|Sc\|^2_2.$$
Therefore, if $\det M\ne 0$, then the operator
$$(UT)^{-1}=\begin{pmatrix} M^{-1} & 0\\-VM^{-1} & 1 \end{pmatrix}$$
is bounded. Hence, $(TU)^{-1} = U^{-1} (UT)^{-1} U$ is bounded and the system $\{e_{\xi_l}\}_{l=1}^m \cup \{e_l\}_{l\in\Z\setminus J}$ is a Riesz basis. To prove the converse, observe that the Riesz basis property implies invertibility of the synthesis operator $T$. But invertibility of $T$ implies surjectivity of $M$, which is equivalent to $\det M \ne 0$ for finite-dimensional operators. \par

We now prove the equivalence of (1) and (3). If some $\xi_j \in \Z$ does not belong to $J$, the corresponding column of $M$ is zero, since $M_{nk} = \sinc(\xi_k-j_n)$. This proves $(1) \rightarrow (3)$. \par

Now consider $\xi_k \in \Z \cap J$. We have $\xi_k = j_s$ for some $s$ and so its column in $M$ has the form $M_{nk} = \sinc(j_s-j_n) = \delta_{n,s}$.  Expand the determinant of $M$ along these columns. Since all such $\xi_k$ are distinct, we have $\det M = \pm \det B$, where $B$ contains no such columns. Note that if $\{\xi_k\} \subset \Z \cap J$, the original determinant equals $1$ or $-1$ and the work is done. Otherwise, after reindexing the remaining rows and columns we have
$$B_{nk} = \sinc (\xi_k-j_n)=\frac{(-1)^{j_n}\sin(\pi\xi_k)}{\pi(\xi_k-j_n)},$$
 where the $\xi_k$ are distinct nonintegers and the $j_n$ are distinct integers. After extracting $\sin(\pi \xi_k) \ne 0$ from the $k$-th column we obtain $\det B = \const \det B'$, where $B'$ is the Cauchy matrix $(1/(\xi_k-j_n))$. Hence, $\det M \ne 0$. 
\end{proof}

Now assume $\det M \ne 0$ and put $\Gram=S^*S$. Consider the system $\{e_{\xi_j}\}_{j = 1}^m \cup \{e_l\}_{l\in\mathbb Z\setminus J}$, where $J=\{j_1,\ldots,j_m\}$. Let $D_l$ be the diagonal projection onto the first $l$ coordinates of $\CC^m$. Fix $C \geq 1$.

\begin{lemma}\label{lem:algprefixes}
Let $P_l$ denote the coordinate projection onto the first $l$ insertions of the system. Then $P_l=SD_l M^{-1}F$ and we have
$$\|P_l\|^2=\frac1{2\pi}\|M^{-*}D_l \Gram D_l M^{-1}\|.$$
In particular, $\|P_l\| \leq C$ if and only if $D_l \Gram D_l \preceq 2\pi C^2M^*M$.
\end{lemma}
\begin{proof}
Since $\det M \ne 0$, the system $\{e_{\xi_l}\}_{l=1}^m \cup \{e_l\}_{l\in\Z \setminus J}$ is a Riesz basis. Now consider 
$$f =  \sum_{j=1}^m c_j e_{\xi_j} + \sum_{l \not\in J} d_l e_l = S(c) + \sum_{l \not\in J} d_l e^{ilt} \in L^2(-\pi,\pi).$$
By the definition of $M$, we have $Mc = F S(c)$. By orthogonality we obtain
$$c = M^{-1} FS(c) = M^{-1}F(f - \sum_{l \not\in J} d_l e^{ilt}) = M^{-1} F(f)$$
and $d_l = \widehat{f}(l) - \sum_{j=1}^m c_j \sinc(\xi_j - l)$ for $l \not \in J$. This proves the formula for $P_l$. Next we have
$$F^*c=\frac1{2\pi}\sum_{l=1}^m c_l e_{j_l} \quad \text{and} \quad FF^*=\frac1{2\pi} 1_m.$$
Put $T_l=SD_l M^{-1}$ and observe that $\|T_l F\|^2= \|T_l F F^* T_l^*\| =\frac1{2\pi}\|T_l^*T_l\|=\frac1{2\pi}\|T_l\|^2$. This proves the norm formula. 
\end{proof}

Let $E\subset \Z \setminus J$ be a finite set of retained frequencies. Order the basis by placing the inserted exponentials first (in the order inherited from $J$), followed by the retained exponentials in any order. Denote by $P_{\mathrm{ins}, E}$ the coordinate projection onto all inserted exponentials and the retained exponentials with frequencies in $E$. Then
$$P_{\mathrm{ins}, E}=Q_{E}+(1-Q_{E})SM^{-1}F.$$
It acts as the identity on $H_E$, as $(1-Q_{E})SM^{-1}F$ on $H_J$, and as zero on $H_{J \cup E}^\perp$. The images of the first two subspaces are orthogonal. Therefore, if $E\ne\varnothing$, we have
$$\|P_{\mathrm{ins},E }\|=\max\left\{1,\frac1{\sqrt{2\pi}}\|(1-Q_{E})SM^{-1}\|\right\}\leq  \max\left\{1, {{1} \over {\sqrt{2\pi}}} \| SM^{-1}\|\right\} = \max\{1, \|P_m\|\} =  \|P_m\|.$$
If $E=\varnothing$, the projection coincides with $P_m$. This proves the following corollary.

\begin{corollary}\label{cor:norm} If the assumptions of Lemma \ref{lem:algprefixes} are valid and $D_l \Gram D_l \preceq 2\pi C^2 M^* M$ for every $0 \leq l \leq m$, then the basis ordered as $e_{\xi_1},\ldots,e_{\xi_m}$ followed by all retained standard vectors in any order has basis constant at most $C$.
\end{corollary}

\subsection{Steps of the construction}

Fix $C\ge 1$ and fix a positive sequence $\eps_n \to 0$. For every $k \in \N$ choose $m_k\ge1$ and choose some integers $j_{k,1}<\ldots<j_{k,m_k}$. Choose some distinct real numbers $\xi_{k,1},\ldots,\xi_{k,m_k}$. The index determines the order of (the future) insertions. \par

{\bf Step 1: matrices.} For each $k$ define the synthesis operator and the coefficient operator
$$S_{k}c=\sum_{l=1}^{m_k} c_l e_{\xi_{k,l}} \quad \text{and } \quad F_{k}f=\frac1{2\pi}\bigl(\langle f,e_{j_{k,l}}\rangle\bigr)_{l=1}^{m_k}.$$
Compute the substitution matrix $M_k$ and the Gram matrix at stage $k$:
$$M_k=F_{k}S_{k}=\bigl(\sinc(\xi_{k,l}-j_{k,d})\bigr)_{d,l=1}^{m_k} \quad \text{and} \quad \Gram_k=S_{k}^*S_{k}=\bigl(2\pi\sinc(\xi_{k,l}-\xi_{k,d})\bigr)_{d,l=1}^{m_k}.$$
These matrices are defined in coordinates related to the centre of the block to be inserted. We have $\|S_{k}c\|^2=c^*\Gram_k c$. For the construction we will need the condition 
$$\det M_k \ne 0.$$ 
By Lemma~\ref{lem:finrep}, each replacement made in the standard basis gives a Riesz basis. Step 8 proves the same for the first $k$ replacements taken together.

{\bf Step 2: projections.} For $0\leq l \leq m_k$ denote by $D_l$ the diagonal projection onto the first $l$ coordinates and set $D_0=0$. We need the following condition:
$$D_l \Gram_k D_l \preceq 2\pi C^2M_k^*M_k\quad\text{ for every }\quad 0\leq l \leq m_k.$$

Lemma~\ref{lem:algprefixes} now gives $\| S_k D_l M_k^{-1} F_k \| \leq C$ for every $0 \leq l \leq m_k$. By Corollary~\ref{cor:norm} this will lead to the Schauder property.

{\bf Step 3: coefficients.} To construct a system with no positive lower Riesz bound we need the condition
$$\lambda_{\min}(\Gram_k)\leq \eps_k^2,$$
where $\lambda_{\min}$ denotes the smallest eigenvalue of a Hermitian matrix. If this condition is satisfied, choose a vector $c_k \in\CC^{m_k}$ from the unit sphere such that $c_k^*\Gram_k c_k \leq \eps_k^2$: for example, one may take a unit eigenvector corresponding to $\lambda_{\min}(\Gram_k)$. The invertibility of $M_k$ implies that $S_{k}$ is injective. Hence, $0<\|S_{k}c_k\|_{2}\leq \eps_k$. \par

{\bf Step 4: radii.}  Set $\Pi_T=Q_{[-T,T]\cap\Z}$. Consider integers $T$ in increasing order, starting with the smallest one satisfying
$$T>\max_{1\leq l\leq m_k} \{|j_{k,l}|,|\xi_{k,l}|\}.$$
For each such $T$ consider the matrix $B_{k,T}=\bigl(\sinc(\xi_{k,n}-l)\bigr)_{\substack{-T\leq l \leq T, 1\leq n\leq m_k}}$. Define $T_k$ to be the first value of $T$ for which
$$M_k^{-*}\bigl(\Gram_k-2\pi B_{k,T}^*B_{k,T}\bigr)M_k^{-1}\preceq 2^{-2k-2} 1_{m_k},$$
where $1_{m_k}$ is the identity matrix on $\CC^{m_k}$. Such a $T_k$ exists: indeed, by Parseval's identity we have
$$\Gram_k-2\pi B_{k,T}^*B_{k,T}=S_{k}^*(1-\Pi_T)S_{k}.$$
Thus, the condition on matrices is equivalent to $\|(1-\Pi_{T_k})S_{k}M_k^{-1}\|\leq 2^{-k-1}$. Now we can find such a $T_k$ for each fixed $k$, since
$$\|(1-\Pi_T)S_{k}M_k^{-1}\|^2\leq \|M_k^{-1}\|^2\sum_{j=1}^{m_k}\|(1-\Pi_T)e_{\xi_{k,j}}\|_2^2$$
and each term in this finite sum tends to zero as $T\to\infty$. Therefore, $\|(1-\Pi_T)S_{k}M_k^{-1}\|\longrightarrow 0$ as $T\to\infty$. The radius is chosen after computing $M_k^{-1}$, since the required tail estimate must account for coefficient recovery by this matrix.

{\bf Step 5: the intervals.} For the first block we put $r_1=T_1+1$ and then by induction
$$r_k=r_{k-1}+T_{k-1}+T_k+1\quad\text{ for }\quad k \geq 2.$$
Define the interval $I_k=[r_k-T_k,r_k+T_k]\cap\Z$ and observe that $\min I_1=1$ and $\min I_{k+1}=\max I_k+1$. In particular, we have
$$\N=\bigcup_{k\ge 1}I_k$$ 
as a disjoint union. The center of $I_k$ is $r_k=\frac{\min I_k+\max I_k}{2}$, and we have $\#I_k=2T_k+1$. The entire family of blocks $\{I_n\}_{n \ge1}$ is now fixed.

{\bf Step 6: removed and inserted frequencies.} Define the set of frequencies to be removed at stage $k$:
$$
J_k=\{r_k+j_{k,n} \mid 1 \leq n\leq m_k\}.
$$
Specify the inserted frequencies using the same labels:
$$
\lambda_{k,r_k+j_{k,n}}=r_k+\xi_{k,n}\quad\text{ for all }\quad 1\leq n\leq m_k.
$$
Let $\Lambda_k=\{\lambda_{k,n} \mid n\in J_k\}$ be the set of inserted frequencies. The order in $\Lambda_k$ is determined by increasing indices $n$. The choice of $T_k$ ensures that $J_k\subset I_k$ and that $\Lambda_k \subset (\min I_k,\max I_k)$. Note that $\#J_k=\#\Lambda_k=m_k$. For the translated block, define
$$S_{k,1}a=\sum_{j\in J_k}a_j e_{\lambda_{k,j}} \quad \text{and} \quad F_{k,1}f=\frac1{2\pi}\bigl(\langle f,e_j\rangle\bigr)_{j\in J_k}.$$
Let $Q_k=Q_{I_k}$, where $I_k$ is defined in Step 5. Identify coordinate $l$ with the index $n=r_k+j_{k,l}$ for $1 \leq l \leq m_k$. For the unitary modulation $(U_rf)(t)=e^{irt}f(t)$ we have $S_{k,1}=U_{r_k}S_{k}$ and $F_{k,1}=F_{k}U_{-r_k}$. Also we have $Q_k=U_{r_k}\Pi_{T_k}U_{-r_k}$. Thus, translation leaves the matrices unchanged:
$$F_{k,1}S_{k,1}=M_k\quad \text{and} \quad S_{k,1}^*S_{k,1}=\Gram_k.$$
The tail condition becomes $\|(1-Q_k)S_{k,1}M_k^{-1}\|\leq 2^{-k-1}$. If Step 3 has been performed, relabel the entries of the vector $c_k$ defined in Step 3 according to
$c_{k,r_k+j_{k,n}}=(c_k)_n$ for all $1\leq n \leq m_k$. Then for this system we have $\sum_{n\in J_k}|c_{k,n}|^2=1$ and $ \left\|\sum_{n\in J_k}c_{k,n}e_{\lambda_{k,n}}\right\|_2\leq \varepsilon_k$.

{\bf Step 7: the finite stage.}  At the $k$-th stage of the algorithm consider the block $I_k$, remove the vectors with frequencies in $J_k$ and insert the vectors with frequencies in $\Lambda_k$. Retain all other frequencies. After the current block has been constructed, stage $k$ contains the replacements in blocks $1,\ldots,k$:
$$\Lambda^{(k)}=\left(\Z\setminus\bigcup_{n=1}^{k}J_n\right)\cup\bigcup_{n=1}^{k}\Lambda_n.$$
At stage $k$, use the following order. For every $n \geq 1$ first list the frequency $1-n$. If $n\leq k$, follow it by $\Lambda_n$ in the order of the labels and then $I_n\setminus J_n$ in increasing order. If $n>k$, follow it by the entire interval $I_n$ in increasing order. Thus the Fourier tail after the first $k$ modified blocks is
$$-k,\ I_{k+1},\ -k-1,\ I_{k+2},\ -k-2,\ I_{k+3},\ldots.$$
The length of the fixed initial segment is $d_k =\sum_{n=1}^{k}(1+\#I_n)=k+r_k+T_k \to \infty$. The first $d_k$ entries remain unchanged at every later stage. Passing from stage $k$ to stage $k+1$ modifies only the list corresponding to $I_{k+1}$. \par
The block introduced at the $k$-th stage is $X_k=\bigl(\{e_{\lambda_{k,n}}\}_{n\in J_k},\{e_n\}_{n\in I_k\setminus J_k}\bigr)$. The ordered system at stage $k$ is
\begin{equation}\label{eq:kbasis}
e_0, \ X_1 , \ e_{-1}, \ X_2, \ \ldots, e_{1-k}, \ X_k, \ e_{-k}, \ \{e_s\}_{s \in I_{k+1}}, \ e_{-k-1},\ \{e_s\}_{s \in I_{k+2}}, \ldots 
\end{equation}

{\bf Step 8: every stage is Riesz.} We now verify the Riesz basis property of the entire system at a finite stage. Set $\mathcal J_k=\bigcup_{n=1}^{k}J_n$. The matrix of the combined replacement is
$$\mathcal M_k=\bigl(\sinc(\lambda_{n,l}-t)\bigr)_{\substack{t\in\mathcal J_k\\1\leq n\leq k,\ l\in J_n}} = \begin{pmatrix}\boxed{M_1} & \star & \star & \cdots & \star \\ \star & \boxed{M_2} & \star & \cdots & \star \\ \star & \star & \boxed{M_3} & \cdots & \star \\ \vdots & \vdots & \vdots & \ddots & \vdots \\ \star & \star & \star & \cdots & \boxed{M_k}\end{pmatrix}.$$
The rows and columns are grouped by blocks. Its diagonal blocks are $M_1,\ldots, M_k$ from Step 1. The off-diagonal blocks are generally nonzero. By Lemma~\ref{lem:finrep}, the system at stage $k$ forms a Riesz basis if and only if $\det\mathcal M_k\ne 0$.\par
It remains to note that the inserted frequencies within each block are distinct. Their placement inside pairwise disjoint intervals excludes coincidences between inserted frequencies from different blocks. The invertibility of $M_n$ for each $n \leq k$ gives $\Lambda_n\cap\Z\subset J_n$ by Lemma~\ref{lem:finrep}. Consequently,
$$\left(\bigcup_{n=1}^{k}\Lambda_n\right)\cap\Z\subset\mathcal J_k.$$
Applying the same lemma to the combined replacement gives $\det \mathcal M_k \ne 0$. Thus the system \eqref{eq:kbasis} is a Riesz basis.

{\bf Step 9: the final system.} The initial segments fixed in Step~7 define the final order:
$$0,\ \Lambda_1, \ I_1\setminus J_1,\ -1,\ \Lambda_2, \ I_2\setminus J_2, \ -2, \ \Lambda_3, \ I_3\setminus J_3, \ -3, \ldots$$
Each finite set is listed in the previously specified order. Denote the corresponding sequence of exponentials by $\{x_j\}_{j\ge1}$. Its frequency set is
$$\Lambda=\left(\Z\setminus\bigcup_{n\ge1}J_n\right)\cup\bigcup_{n\ge 1}\Lambda_n.$$
We first fix the local frequencies and compute the matrices, then choose the radii $T_k$ to control the tails after multiplication by $M_k^{-1}$, and finally place the blocks at the centres $r_k$.

\subsection{Correctness of the Algorithm} Fix $C \ge 1$ and choose a positive sequence $\eps_n \to 0$. Suppose that for every $k\in \N$ the data at the $k$-th stage satisfy
\begin{equation}\label{eq:data}
 \det M_k\ne 0, \quad \lambda_{\min}(\Gram_k)\leq\varepsilon_k^2, \quad D_l\Gram_k D_l \preceq2\pi C^2M_k^*M_k
\end{equation}
for every $0\leq l\leq m_k$ with $C$ and $\eps_k$ fixed above. Let $\{x_j\}_{j=1}^{\infty} \subset L^2(-\pi,\pi)$ be the final system from Step~9. Recall that $Q_n := Q_{I_n} \colon f \mapsto \sum_{j \in I_n} \widehat f(j) e_j$ is the orthogonal projection onto $H_{I_n}$. For each $k\in \N$ consider the finite sequence
\begin{equation}\label{eq:Ybasis}
Y_k=\bigl(\{Q_k e_{\lambda_{k,n}}\}_{n\in J_k},\{e_n\}_{n\in I_k \setminus J_k}\bigr) \subset H_{I_k},
\end{equation}
with both groups listed in their prescribed orders.
\begin{lemma}\label{lem:Ybasis}
The system $Y_k$ forms a basis of $H_{I_k}$ and every initial coordinate projection of $Y_k$ has norm at most $C$.
\end{lemma}
\begin{proof}
 Observe that $\#Y_k = \#I_k$. Order the Fourier coordinates of $H_{I_k}$ as follows: $J_k$ first and $I_k\setminus J_k$ second, in increasing order within each group. Now the coordinate matrix of $Y_k$ in the standard Fourier basis of $H_{I_k}$ is
$$\begin{pmatrix} M_k & 0\\
\star & 1_{\# I_k - m_k} \end{pmatrix}.$$
Since $J_k \subset I_k$, we have $F_{k,1}Q_k=F_{k,1}$ and so $F_{k,1} Q_k S_{k,1}=M_k$. Now $Y_k$ is a basis of $H_{I_k}$, since $F_{k,1}Q_k S_{k,1}=M_k$ is invertible on $\CC^{m_k}$. \par
In order to show that every initial coordinate projection of $Y_k$ has norm at most $C$, consider the basis $B_k$ of $L^2(-\pi,\pi)$ obtained by performing a replacement only in $I_k$, ordered as
$$B_k = \bigl(\{e_{\lambda_{k,n}}\}_{n \in J_k}, \{e_j\}_{j \in I_k\setminus J_k}, \{e_j\}_{j\notin I_k}\bigr),$$
where the first two groups have their prescribed orders and the last group is enumerated arbitrarily. Its initial coordinate projections have norm at most $C$ by Corollary~\ref{cor:norm}. Fix $f\in H_{I_k}$ and write its expansion in this Riesz basis as
$$f=\sum_{n\in J_k}a_n e_{\lambda_{k,n}}+\sum_{j\notin J_k}b_j e_j.$$
Applying $Q_k$ gives
$$f=\sum_{n\in J_k} a_n Q_k e_{\lambda_{k,n}}+\sum_{j\in I_k \setminus J_k}b_j e_j.$$
Since $Y_k$ is a basis, this is the unique expansion of $f$. Let $\Pi_{k,r}$ be the coordinate projection onto the span of the first $r$ vectors of $Y_k$. Let $P_{k,r}^{loc}$ be the initial projection onto the first $r$ vectors of $B_k$. Now we obtain
$$\Pi_{k,r} f= Q_k P_{k,r}^{loc} f \quad \text{ for all } \quad f\in H_{I_k},\ 0\leq r \leq\#I_k.$$
Consequently, $\|\Pi_{k,r}f\|_2 \leq \|P_{k,r}^{loc} f\|_2\leq C\|f\|_2$ and the claim is proved.
\end{proof}

Now we can prove the correctness of the algorithm.

\begin{theorem} \label{th:algorithm}
Let stages of the algorithm satisfy \eqref{eq:data}. Then the final system $\{x_j\}_{j=1}^{\infty}$ is a Schauder basis of $L^2(-\pi,\pi)$ with basis constant at most $3C$ and no positive lower Riesz bound.
\end{theorem}

\begin{proof} The radii are defined by Step 4. For every $k \in \N$ consider the basis $Y_k$ of $H_{I_k}$ defined in \eqref{eq:Ybasis}. Let $\{y_j\}_{j=1}^{\infty}$ be the sequence
$$e_0,\ Y_1,\ e_{-1},\ Y_2,\ e_{-2},\  Y_3,\ldots.$$
The spaces $H_{I_n}$ and $\Span\{e_{1-n}\}$ are pairwise orthogonal, and their algebraic sum is dense in $L^2(-\pi,\pi)$. Recall that $d_k = \sum_{j=1}^k (1+\# I_j)$ and $d_0=0$ (see Step~7). For any $f\in L^2(-\pi,\pi)$ define $\widetilde P_N(f)$ to be the initial projection of $f$ onto $y_1,\ldots,y_N$: if $d_{k-1}< N \leq d_k$ and $r = N - d_{k-1}-1$, then the first $N$ vectors of $\{y_j\}_{j=1}^{\infty}$ consist of $e_0, \ldots, e_{1-k}$, the entire blocks $Y_1, \ldots, Y_{k-1}$ and the first $r$ vectors of $Y_k$. We have
$$\widetilde P_N \colon f \mapsto \sum_{l=1}^k \widehat{f}(1-l) e_{1-l} + \sum_{l=1}^{k-1} Q_l f + \Pi_{k,r} Q_k f.$$
Orthogonality gives 
$$\|\widetilde P_N\|\leq \max\{1, \| \Pi_{k,r} \|\} \leq C.$$
We have $\widetilde P_N \to 1$ strongly, since these operators eventually fix every vector in the algebraic sum. Uniqueness of the expansion follows from uniqueness in each summand. Thus $\{y_j\}_j$ is a Schauder basis with basis constant at most $C$.\par

We now prove that $\{x_j\}_{j=1}^{\infty}$ is a Schauder basis with basis constant at most $3C$. Define
$$T=1+\sum_{n=1}^{\infty} (1-Q_n) S_{n,1} M_n^{-1} F_{n,1}.$$
Using Step 6 and $\|F_{n,1}\|=(2\pi)^{-1/2}$, we obtain that this series converges in operator norm and
$$\|T-1\|\leq \frac1{\sqrt{2\pi}}\sum_{n=1}^{\infty} 2^{-n-1}=\frac1{2\sqrt{2\pi}}<\frac12.$$
Consequently, $T$ is invertible and $\|T\| \cdot \|T^{-1}\| \leq {{1+ \|T-1\|} \over {1-\|T-1\|}} < 3$. For $l\ne k$ we have $F_{l,1} Q_k=0$ and $F_{k,1} Q_k = F_{k,1}$, while $F_{k,1}Q_k S_{k,1}=M_k$ as in Lemma~\ref{lem:Ybasis}. Therefore
$$TQ_k S_{k,1} = Q_k S_{k,1}+ (1-Q_k) S_{k,1} M_k^{-1} F_{k,1} Q_k S_{k,1}= Q_k S_{k,1}+ (1-Q_k) S_{k,1} M_k^{-1} M_k = S_{k,1}.$$
Hence,
$$T Q_k e_{\lambda_{k,n}}=e_{\lambda_{k,n}}\quad \text{ for all } \quad  n\in J_k \text{ and } k \in \N.$$
Moreover, $T$ fixes every $e_m$ with $m\not\in \bigcup_{k \geq 1} J_k$: indeed, in this case we have $F_{k,1} e_m = {{1} \over {2 \pi}} (\langle e_m, e_j \rangle )_{j \in J_k} = (\delta_{m,j})_{j \in J_k} = 0$. Now $x_j=T(y_j)$, so $\{x_j\}_{j=1}^{\infty}$ is a Schauder basis. Let $P_N$ be its $N$-th partial sum projection. Then $P_N=T\widetilde P_N T^{-1}$ and $\|P_N\|\leq 3C$.\par

Finally, take the unit coefficient vectors $\{c_{k,n}\}$ chosen in Step 3 and relabeled in Step 6. Put $v_k= \sum_{n \in J_k} c_{k,n}  e_{\lambda_{k,n}}$. Then $\sum_{n \in J_k} |c_{k,n}|^2=1$ and $0<\|v_k\|_2\leq \eps_k$. These finite combinations occur in the final sequence, so $\eps_k \to 0$ rules out a positive lower Riesz bound.\par
\end{proof}

\section{Existence of the basis} \label{sect:Example}

In this section we construct the example of the basis. We first explain the choice, then give the construction. The required estimates are proved in Section~\ref{sect:Tech}.

\subsection{The motivation}
Suppose we have already performed all replacements and obtained the final system. To ensure that this system is not a Riesz basis, it is sufficient to show that $\lambda_{\min}(\Gram_k) \to 0$. In order to obtain the Schauder basis property, we need uniformly bounded initial projections for each local replacement. By Lemma~\ref{lem:algprefixes} we need an absolute constant $C_1$ such that $\|S_k D_l M_k^{-1} F_k \| \leq C_1$ for every $k$ and $0 \leq l \leq m_k$. Taking $l = m_k$ gives $\|S_k M_k^{-1}\| \leq \sqrt{2\pi} C_1$. Since
$$\lambda_{\min} (\Gram_k) \|M_k^{-1}\|^2 \leq \|S_k M_k^{-1}\|^2 \leq 2 \pi C_1^2,$$
it is enough to make $\|M^{-1}_k\| \to \infty$ while keeping all initial projection norms uniformly bounded. \par

Consider one finite replacement as in Subsection~\ref{sect:finrep}. Fix $m \in \N$, let $j_1< \ldots < j_m$ be integer frequencies to be replaced, and let $\xi_1, \ldots, \xi_m \in \R$ be distinct frequencies to be inserted. Let $S$, $M$ and $\Gram = S^* S$ be the corresponding operators defined in Subsection \ref{sect:finrep}, and assume $\det M \ne 0$. The generating function of this replacement is
$$G(z) = {{\sin \pi z} \over {\pi}} \prod_{l=1}^{m} {{z-\xi_l} \over {z-j_l}}.$$
A partial fraction decomposition gives
$${{G(z)} \over {G'(\xi_l) (z-\xi_l)}} = \sum_{n=1}^m (M^{-1})_{ln} \sinc(z-j_n).$$
The left hand side corresponds to the biorthogonal vector associated with $e^{i \xi_l t}$ in the replacement basis. Let $R \geq 1$ and choose $E \subset \{1, \ldots, m\}$ such that $|\xi_l| \leq R$ for every $l \in E$. Then we have
$$|G(iR)|^2 \sum_{l \in E} {{1} \over {|G'(\xi_l)|^2 |iR-\xi_l|^2}} \leq \|M^{-1}\|^2 \sum_{n=1}^m |\sinc(iR-j_n)|^2 \leq \|M^{-1}\|^2 {{\|e_{iR}\|_2^2} \over {2\pi}} \leq  \|M^{-1}\|^2 {{e^{2\pi R}} \over {4\pi R}}.$$
Since $|iR-\xi_l|^2\leq 2R^2$, we have the desired lower bound for $\|M^{-1}\|$:
$$\|M^{-1}\|^2 \geq {{e^{-2\pi R} |G(iR)|^2} \over {R}} \sum_{l \in E} {{1} \over {|G'(\xi_l)|^2}}.$$
For the $k$-th block, put $R = 3^{k+1}$ and choose $E_k \subset \{1, \ldots, m_k\}$ with $|\xi_{k,l}| \leq R$ for $l \in E_k$. We seek $|G_k(iR)| \asymp e^{\pi R}$ as $R \to \infty$, with constants independent of $k$. We aim for
\begin{equation}\label{eq:desired}
 |G_k'(\xi_{k,l})|  \lesssim R^{-\alpha} \quad \text{and} \quad \#E_k \asymp R^d \quad \text{as } R \to \infty
\end{equation}
with constants independent of $k$ and $l \in E_k$. In this case 
$$\|M_k^{-1}\|^2 \gtrsim  R^{d+2\alpha-1}.$$
This gives $\lambda_{\min}(\Gram_k) \to 0$ whenever $d+2\alpha > 1$. We take $\#E_k = 2^{k+1}$ and $\alpha = 1/4$.  Thus $d = \log_3 2$ and $d + 2\alpha >1$. Then
$$\|M_k^{-1}\|^2 \gtrsim {{2^{k+1}} \over {\sqrt{R} }} = \left( {{2} \over {\sqrt{3} }} \right)^{k+1} \to \infty.$$
Later we prove that $|G_k'(\xi_{k,l})| \leq |G_k(\xi_{k,l} +i)|$, so we need small values of $|G_k|$ at height 1 near many frequencies, while keeping large values of $|G_k|$ at height $R$. Furthermore,
$$\ln \left| {{G(z)} \over {\sin \pi z/\pi}} \right| = \sum_{l=1}^m \ln |z-\xi_l | - \sum_{l=1}^m \ln |z-j_l| = \int_{\R} \ln |z-t| \ d\mu(t),$$
where $\mu = \sum_{l=1}^m \delta_{\xi_l} - \sum_{l=1}^m  \delta_{j_l}$. In order to find such a function $G$, consider the continuous version of the problem:
$$u(x,y) = \int_{\R} \ln |x+iy - t| g(t) \ dt.$$

Since $\mu(\R) = 0$, we must impose the condition $\int_{\R} g(t) \ dt =0$. We prescribe a continuous compactly supported boundary function $\varphi$ and seek its harmonic extension, since the logarithmic kernel is harmonic in the upper half-plane:
\begin{equation}\label{eq:gint}
u(x,y) = P_y \varphi =  {{1} \over {\pi}} \int_{\R} {{y \varphi(t)} \over {(x-t)^2 + y^2}} \ dt \quad \text{for } y > 0.
\end{equation}
Our boundary function $\varphi$ is continuous, but ${{\partial} \over {\partial y}} P_y \varphi$ may be unbounded as $y \to 0$. We therefore smooth $\varphi$: fix some $a>0$ and consider
$$u_a(x,y) = P_y (P_a \varphi) (x) = P_{y+a} \varphi(x).$$
A standard calculation gives
\begin{equation}\label{eq:gresult}
g_a(x) = {{1} \over {\pi}} {{\partial u_a} \over {\partial y}} (x,0) = {{1} \over {\pi}}  {{\partial} \over {\partial y}} (P_y \varphi) (x) |_{y=a}.
\end{equation}

The function $g_a$ has integral zero and satisfies
$$\int_{\R} \ln |x + i y - t| g_a(t) \ dt = (P_{y+a} \varphi)(x) \quad \text{for } y > 0.$$

The specific boundary function chosen below will give the small derivatives needed above, and its harmonic extension will control the contour integrals used to bound all initial projections in Subsection~\ref{sect:schauder}.

\subsection{The construction} The construction starts with a finite set of integers
\begin{equation}\label{eq:Ak}
 A_k = \left\{\pm\sum_{n=0}^{k} 2 \omega_n 3^n \mid \omega_n \in \{0,1\}\right\}.
\end{equation}

We schematically draw the first three such sets:
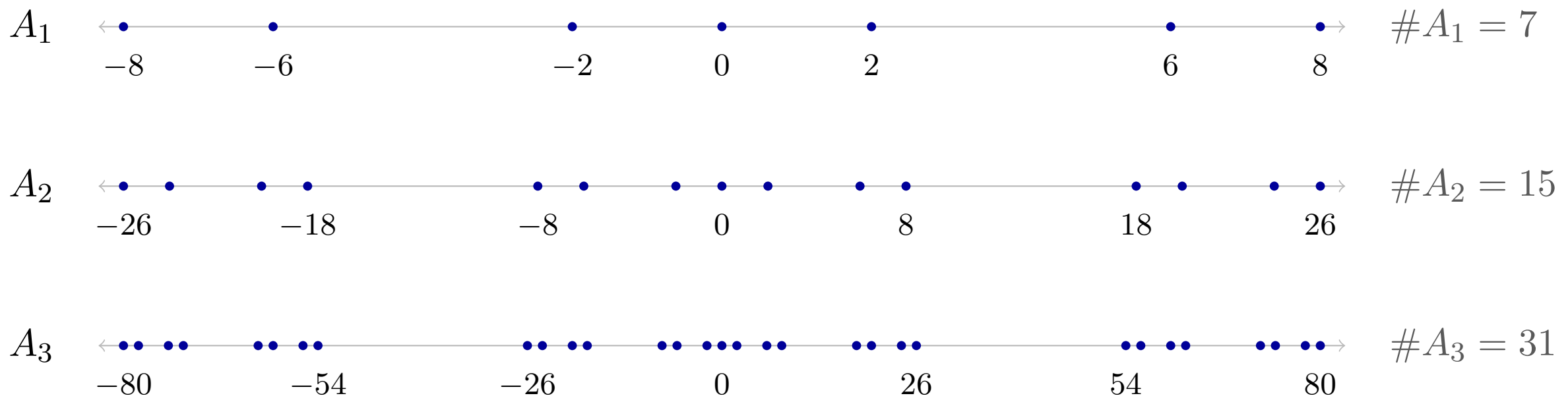
\begin{figure}[htbp]
\centering
\begin{tikzpicture}[x=1cm, y=1cm, font=\small]

\foreach \k/\y/\n in {1/3.2/7, 2/1.6/15, 3/0/31} {
    \draw[gray!55, <->, line width=0.4pt] (-6.25, \y) -- (6.25, \y);
    \node[anchor=east] at (-6.55, \y) {$A_{\k}$};
    \node[anchor=west, text=gray!70!black] at (6.55, \y) {$\# A_{\k} = \n$};
}

\foreach \x in {-8, -6, -2, 0, 2, 6, 8} {
    \fill[blue!60!black] ({6 * \x / 8}, 3.2) circle (1.3pt);
    \node[below=4pt, font=\scriptsize] at ({6 * \x / 8}, 3.2) {$\x$};
}

\foreach \x in {-26, -24, -20, -18, -8, -6, -2, 0, 2, 6, 8, 18, 20, 24, 26} {
    \fill[blue!60!black] ({6 * \x / 26}, 1.6) circle (1.3pt);
}
\foreach \x in {-26, -18, -8, 0, 8, 18, 26} {
    \node[below=4pt, font=\scriptsize] at ({6 * \x / 26}, 1.6) {$\x$};
}

\foreach \x in {-80, -78, -74, -72, -62, -60, -56, -54, -26, -24, -20, -18, -8, -6, -2, 0, 2, 6, 8, 18, 20, 24, 26, 54, 56, 60, 62, 72, 74, 78, 80} {
    \fill[blue!60!black] ({6 * \x / 80}, 0) circle (1.3pt);
}
\foreach \x in {-80, -54, -26, 0, 26, 54, 80} {
    \node[below=4pt, font=\scriptsize] at ({6 * \x / 80}, 0) {$\x$};
}

\end{tikzpicture}
\caption{The sets $A_1$, $A_2$, and $A_3$}
\label{fig:Ak}
\end{figure}

The set $A_k$ is symmetric and its nonnegative part $A_k \cap [0,\infty)$ contains $2^{k+1}$ points in $[0,3^{k+1}-1]$. They are separated by at least 2. We associate with this set the function
\begin{equation}\label{eq:phi}
\varphi_k(x) = \frac14 \ln \min\left\{1, {{1+\dist(x, A_k)} \over {3^{k+1}}}\right\}.
\end{equation}
This is our choice of the boundary function $\varphi$ from the motivation above. At each point $x\in A_k$ we have $\varphi_k(x) = -{{1} \over {4}} \ln 3^{k+1}$, corresponding to the choices $R=3^{k+1}$ and $\alpha = 1/4$ in the motivation. It is chosen to make $|G_k|$ small at height 1 near $A_k$, while keeping $|G_k(i3^{k+1})| \asymp e^{\pi 3^{k+1}}$. The ternary structure serves to control the Poisson averages. This function is nonpositive and vanishes outside $[-2\cdot 3^{k+1}, 2\cdot 3^{k+1}]$. Recall that the Poisson integral of $f$ is defined by
$$(P_y f)(x) = \frac1\pi \int_{\mathbb R} \frac{y f(t)}{(x-t)^2+y^2}\,dt.$$
As in \eqref{eq:gresult}, set
$$g_k(x) = \left. \frac1\pi \frac{\partial}{\partial y}(P_y \varphi_k)(x)\right|_{y=4}.$$
The value $y=4$ is chosen only to make the calculation simple. We now define the $k$-th density $V_k$ as follows: the function $V_k-1$ is obtained by truncating $g_k$ to $[-4 \cdot 3^{k+1}, 4 \cdot 3^{k+1}]$ and subtracting its mean on the interval:
\begin{equation} \label{eq:density}
V_k(x) -1  = \begin{cases} \displaystyle g_k(x)-\frac{1}{8\cdot3^{k+1}}\int_{-4\cdot3^{k+1}}^{4\cdot3^{k+1}} g_k(t) \,dt, & |x| \leq 4\cdot 3^{k+1}, \\
0, & |x| > 4\cdot 3^{k+1}.
\end{cases}
\end{equation}
This choice ensures that $\int_{\R} (V_k(t) - 1) \, dt = 0$, as required in the motivation above. Now let
$$N_k(x) = \int_0^x V_k(t)\,dt$$
be the corresponding counting function. The density is even and satisfies $5/6 < V_k(x) < 7/6$ by Lemma~\ref{lem:density}. Hence, $N_k$ is a strictly increasing continuous bijection of $\mathbb R$, and its inverse $N_k^{-1}$ exists. Since $V_k$ is even and $\int_{\R} (V_k(t)-1) \ dt = 0$, we have $N_k(x) = x$ for $|x| \geq 4\cdot 3^{k+1}$. Hence, $N_k^{-1}(u) = u$ for $|u| \geq 4 \cdot 3^{k+1}$.\par

Set $m_k = 8\cdot 3^{k+1}+1$ and define 
$$j_{k,l} = -4\cdot 3^{k+1}+l-1 \quad \text{ for } 1 \leq l \leq m_k.$$
For each $j \in \Z$ the interval $[N_k^{-1}(j - 1/2), N_k^{-1}(j + 1/2)]$ has measure one with respect to $V_k(t) \, dt$. Now choose the inserted frequencies as the centres of mass of these intervals:
$$\xi_{k,l} = \int_{j_{k,l}-1/2}^{j_{k,l}+1/2} N_k^{-1}(u)\,du = \int_{N_k^{-1}(j_{k,l}-1/2)}^{N_k^{-1}(j_{k,l}+1/2)} t V_k(t) \ dt \quad  \text{for } 1 \leq l \leq m_k.$$
This choice gives a uniform bound on the error when replacing sums of logarithms by integrals in Lemma~\ref{lem:genprofile}, leading to the desired asymptotics for the generating function $G_k$. At the $k$-th step, every inserted frequency lies in $(-4\cdot 3^{k+1}-1, 4\cdot 3^{k+1}+1)$, since $N_k^{-1}$ is strictly increasing and acts as the identity for $|u| \geq 4 \cdot 3^{k+1}$.\par

The frequencies are strictly increasing and satisfy $6/7 \leq \xi_{k,l+1}-\xi_{k,l} \leq 6/5$ for $1 \leq l < m_k$. We use these choices of $m_k$, $j_{k,l}$ and $\xi_{k,l}$ in Step~1, with $M_k$ and $\Gram_k$ defined there. \par

We fix $C = \max\{C_1, C_2\}$, where $C_1, C_2 \geq 1$ are the absolute constants specified in \eqref{eq:Schauder} and \eqref{eq:c2}, and $ \eps_k = C \left(\frac34\right)^{{{k+1} \over {4}}}$ for $k \in \N$.\par

\begin{theorem}\label{th:exist}
For every $k \in\N$ we have $\det M_k \ne 0$ and $\lambda_{\min}(\Gram_k) \leq \varepsilon_k^2$. Moreover,
$$D_l\Gram_k D_l \preceq 2\pi C^2 M_k^* M_k \quad \text{for every } 0 \leq l \leq m_k.$$
\end{theorem}

The proof is given in Section~\ref{sect:Tech}. We now apply Steps~4--9. Choose $T_k$ as in Step 4 and form the intervals $I_k$ with centres $r_k$. As in Step~6, let 
\begin{equation}\label{eq:defJ}
J_k = \{r_k+j_{k,l} \mid 1 \leq l \leq m_k\} \quad \text{and} \quad \lambda_{k,r_k+j_{k,l}} = r_k+\xi_{k,l}. 
\end{equation}

Put $\Lambda_k = \{\lambda_{k, n} \mid n \in J_k\}$.  List each $\Lambda_k$ in increasing order of its labels and each $I_k \setminus J_k$ in increasing order. The final order is
$$0,\ \Lambda_1,\ I_1 \setminus J_1,\ -1,\ \Lambda_2,\ I_2\setminus J_2,\ -2,\ \Lambda_3,\ I_3\setminus J_3,\ldots$$
By Theorem~\ref{th:algorithm}, the resulting exponentials $\{x_j\}_{j=1}^{\infty}$ form a Schauder basis of $L^2(-\pi,\pi)$ with basis constant at most $3C$ and no positive lower Riesz bound.\par

\section{Proof of Theorem~\ref{th:exist}}\label{sect:Tech}

We retain the indexing from Step~6: $\lambda_{k, n}$ replaces the integer frequency $n \in J_k$, and $\lambda_{k, r_k + j_{k, l}} = r_k + \xi_{k, l}$ for $1 \leq l \leq m_k$. We prove the estimates for $\xi_{k, l}$ before translation, since the unitary modulation $U_{r_k}$ leaves the matrices and operator norms unchanged. Hereinafter we fix $k \in \N$. Unless stated otherwise, the constants in asymptotic inequalities are absolute.

\subsection{The density properties}

\begin{lemma}\label{lem:density}
The density defined in \eqref{eq:density} satisfies $5/6 < V_k < 7/6$. We have $N_k(x) = x$ for $|x| \geq 4 \cdot 3^{k + 1}$. The inserted frequencies satisfy $\frac{6}{7} \leq \xi_{k, l + 1} - \xi_{k, l} \leq \frac{6}{5}$ for every $1 \leq l < m_k$. The gaps between the inserted block and the nearest retained integer frequencies satisfy
$$\frac{6}{7} \leq \xi_{k, 1} + 4 \cdot 3^{k + 1} + 1 \leq \frac{6}{5} \quad \text{and} \quad \frac{6}{7} \leq 4 \cdot 3^{k + 1} + 1 - \xi_{k, m_k} \leq \frac{6}{5}.$$
All remaining gaps are of length 1. In particular, $\det M_k \ne 0$.
\end{lemma}

\begin{proof}
Using \eqref{eq:phi} we get $|\varphi_k(x) - \varphi_k(t)| \leq \frac{1}{4} \ln(1 + |x - t|)$: indeed, $1+\dist(t,A_k) \leq (1+\dist(x,A_k))(1+|x-t|)$ by Lipschitz, which gives the claim. We have 
$$g_k(x) = \left. {{1} \over {\pi^2}} {{\partial} \over {\partial y}} \int_{\R} {{y \varphi_k(t)} \over {(x-t)^2+y^2}} dt \right|_{y=4} = {{1} \over {\pi^2}} \int_{\R} \varphi_k(t) {{(x-t)^2-16} \over {((x-t)^2+16)^2}} \ dt.$$
Since the differentiated Poisson kernel has integral zero, we obtain
$$|g_k(x)| = \left| {{1} \over {\pi^2}} \int_{\R} (\varphi_k(t) - \varphi_k(x)) {{(x-t)^2-16} \over {((x-t)^2 +16)^2}} \ dt \right| \leq {{1} \over {\pi^2}} \int_{\R} {{ |\varphi_k(t) - \varphi_k(x)|} \over {(x-t)^2+16}} \ dt \leq $$
$$\leq \frac{1}{8 \pi^2} \int_0^\infty \frac{\ln(1 + 4 s)}{1 + s^2} \, d s \leq \frac{1}{8 \pi^2} \int_0^\infty \frac{\ln 5 (1 + s)}{1 + s^2} \, d s < \frac{1}{12},$$
where in the second inequality we use $s = (t-x)/4$ and evenness. This yields the claimed bounds for $V_k$. Note that $V_k$ is even since $\varphi_k$ is even. The evenness of $V_k$ and the identity $\int_{\mathbb R} (V_k - 1) = 0$ give $N_k(x) = x$ outside $[-4 \cdot 3^{k + 1}, 4 \cdot 3^{k + 1}]$.\par

The density bounds imply $6/7 \leq N_k^{-1}(u + 1) - N_k^{-1}(u) \leq 6/5$: indeed, if $a = N_k^{-1}(u)$ and $b=N_k^{-1}(u+1)$, then $1 = N_k(b) - N_k(a) = \int_a^b V_k(t) \ dt$, which gives $\frac{5}{6}(b-a) < 1<\frac{7}{6}(b-a)$. Now for the inner frequencies we have 
$$\frac{6}{7} \leq \xi_{k,l+1} - \xi_{k,l} = \int_{j_{k,l}-1/2}^{j_{k,l}+1/2} (N_k^{-1}(u+1) - N_k^{-1}(u)) du \leq {{6} \over {5}}.$$
For the boundary frequencies we have 
$$\xi_{k,1} - j_{k,1} +1 = \int_{j_{k,1} - 3/2}^{j_{k,1} - 1/2} (N_k^{-1}(u+1) - N_k^{-1}(u)) du =  \int_{j_{k,1} - 3/2}^{j_{k,1} - 1/2} (N_k^{-1}(u+1) - u) du,$$
$$j_{k,m_k} + 1 - \xi_{k,m_k}= \int_{j_{k,m_k} - 1/2}^{j_{k,m_k} + 1/2}  (N_k^{-1}(u+1) - N_k^{-1}(u)) du = \int_{j_{k,m_k} - 1/2}^{j_{k,m_k} + 1/2}  (u+1 - N_k^{-1}(u)) du,$$
which give the bounds. Since $N_k^{-1}$ is strictly monotone, we have $N_k^{-1} (j_{k,l}-1/2) < \xi_{k,l} < N_k^{-1} (j_{k,l}+1/2)$ by the definition of $\xi_{k,l}$. Now all inserted frequencies lie in $(-4 \cdot 3^{k + 1} - 1/2, 4 \cdot 3^{k + 1} + 1/2)$. All integers in this interval lie in $\{j_{k,1}, \ldots, j_{k, m_k} \}$ by the definition of $j_{k,l}$. Hence, if an inserted frequency is an integer, then it coincides with some $j_{k,s}$. Now Lemma~\ref{lem:finrep} gives $\det M_k \ne 0$.
\end{proof}

\subsection{The Poisson estimates}

We begin with the uniform neighbourhood estimate for the Poisson integral. First we consider the set $A_k$ defined in \eqref{eq:Ak} and define $A_k^+ = A_k \cap [0,+\infty)$. We have $A_k = A_k^+ \cup (-A_k^+)$. 

\begin{lemma}\label{lem:ternary-counting}
Fix an interval $I\subset \R$ of length $r$. For every $0 < h \leq r$ we have
$$|\{x \in I \mid \dist(x, A_k) < h\}| \lesssim r^{\log_3 2} h^{1 - \log_3 2}.$$
\end{lemma}

\begin{proof}
Suppose $I = (a,b)$ and fix $0<h\leq r$. First suppose that $r < 1$. Now $(a-h,b+h)$ contains at most two points of $A_k$, since $A_k \subset 2\Z$ and the length of $(a-h,b+h)$ is $r+2h \leq 3r <3$.  The neighbourhoods of these points have total length at most $4h$. But $4h \leq 4 r^{\log_3 2} h^{1 - \log_3 2}$ in this case. \par

For $r \geq 1$ it is sufficient to prove the estimate only for $A_k^+$.  Choose $q \in \Z$ such that $r \leq 3^q < 3 r$. The interval $(a-h,b+h)$ meets a bounded number of intervals of the form $[j 3^q, (j+1)3^q)$ for $j \in \Z$. Now it is sufficient to consider only points from $A_k^+ \cap [j3^q, (j+1)3^q)$ for each such $j$. \par

If $t = \sum_{n=0}^k 2\omega_n 3^n \in A_k^+ \cap [j3^q, (j+1)3^q)$ then put $\omega_{k+1} = \ldots = \omega_{q-1} = 0$ if $q-1>k$ and observe 
$$t = j3^q + \sum_{n=0}^{q-1} 2\omega_n 3^n.$$
 First consider the case $h \geq 1$. One can find $p \in \{0, \ldots, q\}$ such that $3^{q-p} \leq h < 3^{q-p+1}$. We have 
$$0 \leq t - j 3^q - \sum_{n = q-p}^{q-1} 2\omega_n 3^n = \sum_{n=0}^{q-p-1} 2\omega_n 3^n \leq 3^{q-p}-1.$$
Fix $\beta =  j 3^q +  \sum_{n = q-p}^{q-1} 2\omega_n 3^n$. The preceding inequalities imply $t \in [\beta, \beta + 3^{q-p})$. Each of the numbers $\omega_{q-p}, \ldots , \omega_{q-1}$ has only 2 possible values. Hence, there are at most $2^p$ possible values of $\beta$. Thus every $t$ in the intersection lies in the union of at most $2^p$ intervals of length $3^{q-p} \leq h$. Now the $h$-neighbourhood of every such interval has length at most $3h$ and 
$$\left| \bigcup_{t \in A_k^+ \cap [j 3^q, (j+1)3^q)} (t-h, t+h) \right| \leq 3 \cdot 2^p h.$$
Since $h < 3^{q-p+1}$, we obtain 
$$2^p h = (3^p)^{\log_3 2} h \lesssim r^{\log_3 2} h^{1 - \log_3 2}.$$
The case $0<h<1$ is obvious since in every half-interval of length $3^q$ there are at most $2^q$ points from $A_k^+$, and their $h$-neighbourhoods have total length at most $2\cdot 2^q h \lesssim r^{\log_3 2} h \leq r^{\log_3 2} h^{1-\log_3 2}$.
\end{proof}

\begin{lemma}\label{lem:poisson}
There exists an absolute constant $c_1>0$ such that for every $x \in \mathbb R$ and $y \geq 1$ we have
$$\left|P_{y + 4} \varphi_k(x) - \frac{1}{4} \ln \min\left\{1, \frac{1 + \dist(x, A_k) + y}{3^{k + 1}}\right\} \right| \leq c_1.$$
\end{lemma}

\begin{proof}
Fix $k \geq 1$, $x \in \mathbb R$ and $y \geq 1$, and put $\rho = 1 + \dist(x, A_k) + y$. We first estimate $P_y \varphi_k(x)$. Define
$$\Delta_k(t) = 4 \varphi_k(t) - \ln \min\left\{1, \frac{\rho}{3^{k + 1}}\right\}$$
and put $ \Delta_{k, +}(t) = \max\{\Delta_k(t), 0\}$ and $\Delta_{k, -}(t) = \max\{-\Delta_k(t), 0\}$. We have 
$$\left|4 P_y \varphi_k(x) - \ln \min\left\{1, \frac{\rho}{3^{k + 1}}\right\}\right| \leq  \int_{\mathbb R} \frac{y (\Delta_{k, +}(t) + \Delta_{k, -}(t))}{\pi ((x - t)^2 + y^2)} \ dt.$$
Replacing $y$ by $y + 4$ changes the logarithmic term by at most $\ln 3$. Hence, it is sufficient to show that 
\begin{equation} \label{eq:deltas}
\int_{\R} {{\Delta_{k,+}(t) y \ dt} \over {\pi((x-t)^2+y^2)}} \lesssim 1 \quad \text{and} \quad \int_{\R} {{\Delta_{k,-}(t) y \ dt} \over {\pi((x-t)^2+y^2)}} \lesssim 1.
\end{equation}

Observe that $\dist(t, A_k) \leq \dist(x, A_k)+|t-x|$ for every $x,t\in \R$. Hence, $\Delta_{k,+}(t) \leq \ln \left(1+ {{|t-x|} \over {y}}\right)$. Using $t = x+yu$ we obtain
$${{1} \over {\pi}} \int_{\R} {{\Delta_{k,+}(t) y} \over {(x-t)^2+y^2}} \ dt \leq {{1} \over {\pi}} \int_{\R} {{y} \over {(x-t)^2 + y^2}} \ln \left(1+ {{|t-x|} \over {y}}\right) \ dt = {{1} \over {\pi}} \int_{\R} {{\ln(1+|u|)} \over {1+u^2}} \ du < \infty.$$
It remains to prove the second inequality in \eqref{eq:deltas}.\par

To estimate the negative part, we first establish an estimate over the sets $E_{k, h} = \{t \in \mathbb R \mid 1 + \dist(t, A_k) < h\}$ for $0 < h \leq \rho$. Split the real line: $\R = [x - \rho, x + \rho] \cup \bigcup_j R_j$, where
$$R_j = \{t \in \R \mid 2^j \rho < |t-x| \leq 2^{j+1} \rho\} \quad \text{for } j \geq 0.$$
Each $R_j$ consists of two intervals of length $2^j \rho$. Observe that $E_{k, h} \subset \{t \mid \dist(t, A_k) < h\}$. \par
First assume that $h \leq \rho/4$. In this case one may apply Lemma~\ref{lem:ternary-counting} separately to both intervals that form $R_j$. We obtain
$$|E_{k, h} \cap  R_j| \lesssim (2^j \rho)^{\log_3 2} h^{1 - \log_3 2}.$$
The same lemma applied to $[x-\rho, x+\rho]$ gives $|E_{k, h} \cap [x - \rho, x + \rho]| \lesssim \rho^{\log_3 2} h^{1 - \log_3 2}$. Hence, on $E_{k, h} \cap [x-\rho, x+\rho]$ the kernel satisfies
$${{y} \over {(x-t)^2+y^2}} \lesssim {{1} \over {\rho}},$$
because $\rho \leq |x-t| + h + y$ and so $|x-t| + y \ge 3 \rho/4$. On $E_{k, h} \cap R_j$ the kernel bound is $2^{-2j}/\rho$. Now for $0< h< \rho/4$ we obtain
\begin{equation}\label{eq:ekh}
\int_{E_{k, h}} \frac{y}{\pi ((x - t)^2 + y^2)} \, d t \lesssim \left(\frac{h}{\rho}\right)^{1 - \log_3 2} \left(1 + \sum_{j = 0}^{\infty} 2^{j (\log_3 2 - 2)}\right) \lesssim \left(\frac{h}{\rho}\right)^{1 - \log_3 2}.
\end{equation}
For $\rho/4 < h \leq \rho$, the integral on the left is at most $1 \leq 4^{1 - \log_3 2} (h/\rho)^{1 - \log_3 2}$. Thus \eqref{eq:ekh} holds for every $0 < h \leq \rho$.\par

For every $s > 0$ the definition of $\Delta_k$ gives
$$\{t \mid \Delta_{k, -}(t) > s\} = \left\{t \mid 1 + \dist(t, A_k) < \min\{3^{k + 1}, \rho\} \cdot e^{-s}\right\}.$$
Indeed, $\Delta_{k, -}(t) > s$ is equivalent to $\min\{3^{k + 1}, 1 + \dist(t, A_k)\} < \min\{3^{k + 1}, \rho\}\cdot e^{-s}$. The right-hand side is strictly smaller than $3^{k + 1}$, so the minimum on the left can be removed.\par

Now $h = \min\{3^{k+1}, \rho\} e^{-s} \leq \rho e^{-s}$ for every $s>0$. Applying \eqref{eq:ekh}, we obtain
$$\int_{\mathbb R} \frac{y \Delta_{k, -}(t)}{\pi ((x - t)^2 + y^2)} \, d t = \int_0^\infty \left(\int_{\{t \mid \Delta_{k, -}(t) > s\}} \frac{y}{\pi ((x - t)^2 + y^2)} \, d t\right) ds = $$
$$= \int_0^\infty \left(\int_{E_{k, \min\{3^{k+1}, \rho\} e^{-s}}} \frac{y}{\pi ((x - t)^2 + y^2)} \, d t\right) ds  \lesssim \int_{0}^{\infty} e^{-(1 - \log_3 2) s} \, d s \lesssim 1.$$
\end{proof}

\begin{corollary}\label{cor:phiint} We have $\displaystyle \int_{\mathbb R} |\varphi_k(t)| \, d t \lesssim 3^{k + 1}$.
\end{corollary}
\begin{proof} Apply Lemma~\ref{lem:poisson} for $x = 0$ and $y = 3^{k+1}$. We obtain $|P_{3^{k+1}+4} \varphi_k(0) | \lesssim 1$. Hence,
$${{1} \over {3^{k+1}}} \int_{\R} |\varphi_k(t)| \ dt =\int_{-2 \cdot 3^{k+1}}^{2\cdot 3^{k+1}}{{|\varphi_k(t)| } \over {3^{k+1}}}  \ dt \asymp \int_{-2 \cdot 3^{k+1}}^{2\cdot 3^{k+1}} {{(3^{k+1} + 4)  |\varphi_k(t)| } \over {\pi(t^2+(3^{k+1} + 4)^2)}}\ dt \lesssim - P_{3^{k+1}+4} \varphi_k(0) \lesssim 1.$$

\end{proof}

\subsection{The generating function}

Define $\mathcal Z_k = \bigl(\mathbb Z \setminus \{j_{k, 1}, \ldots, j_{k, m_k}\}\bigr) \cup \{\xi_{k, 1}, \ldots, \xi_{k, m_k}\}$ and let $G_k$ be its generating function, that is
$$G_k(z) = \frac{\sin \pi z}{\pi} \prod_{l = 1}^{m_k} \frac{z - \xi_{k, l}}{z - j_{k, l}}.$$
By Lemma~\ref{lem:density}, the zeros of $G_k$ are precisely the points of $\mathcal Z_k$. These are simple and have consecutive gaps between $6/7$ and $6/5$.

\begin{lemma}\label{lem:genprofile}
For every $x \in \mathbb R$ and $y \geq 1$ we have
$$|G_k(x + i y)| \asymp e^{\pi y} \min\left\{1, \left(\frac{1 + \dist(x, A_k) + y}{3^{k + 1}}\right)^{1/4}\right\}.$$
\end{lemma}

\begin{proof}
It is sufficient to prove $\left| \ln |G_k(x + i y)| - \pi y - P_{y + 4} \varphi_k(x) \right| \lesssim 1$ and apply Lemma~\ref{lem:poisson} to get the claim. For $1 \leq l \leq m_k$ put $a_{k,l} = N_k^{-1}\left(j_{k,l} - {{1} \over {2}}\right)$ and $ b_{k,l} = N_k^{-1}\left(j_{k,l} + {{1} \over {2}}\right)$. Fix $z = x + i y$ for $x \in \R$ and $y \geq 1$. \par

{\bf Step 1}. We first show that 
\begin{equation}\label{eq:Gint}
\begin{aligned}
\sum_{l=1}^{m_k} \bigg|\ln|z-\xi_{k,l}| - \int_{a_{k,l}}^{b_{k,l}} \ln |z-t|  V_k(t) dt \bigg| \lesssim {{1} \over {y}}, \\
\sum_{l=1}^{m_k} \bigg|\ln|z-j_{k,l}| - \int_{j_{k,l}-1/2}^{j_{k,l}+1/2} \ln |z-t| dt \bigg| \lesssim {{1} \over {y}}.
\end{aligned}
\end{equation}

Recall that $\xi_{k,l} \in [a_{k,l}, b_{k,l}]$ and $b_{k,l} - a_{k,l} \leq 6/5$ by the proof of Lemma~\ref{lem:density}. Consider $f_z(t) = \ln |z - t|$. Its second derivative satisfies $|f_z''(t)| \leq |z - t|^{-2}$. If $t, r \in [a_{k,l}, b_{k,l}]$, then $|z - t| \leq |z - r| + 6/5 \lesssim |z - r|$, because $|z - r| \geq y \geq 1$. Taylor's formula therefore yields
$$|f_z(t) - f_z(\xi_{k,l}) - f_z'(\xi_{k,l}) (t - \xi_{k,l})| \lesssim \frac{(t - \xi_{k,l})^2}{|z - t|^2} \lesssim \frac{1}{|z - t|^2}.$$
We obtain
$$\left|\ln |z - \xi_{k,l}| - \int_{a_{k,l}}^{b_{k,l}} \ln |z - t| V_k(t) \ d t\right| = \left| \int_{a_{k,l}}^{b_{k,l}} \left( \ln |z-t| - \ln |z - \xi_{k,l}| - f_z'(\xi_{k,l})(t-\xi_{k,l}) \right) V_k(t) \ dt \right|$$
$$\lesssim \int_{a_{k,l}}^{b_{k,l}} \frac{V_k(t)}{|z - t|^2}  \ d t \leq {{7} \over {6}} \int_{a_{k,l}}^{b_{k,l}} \frac{d t}{|z - t|^2}.$$
The linear term disappears because $\int_{a_{k,l}}^{b_{k,l}} (t - \xi_{k,l}) V_k(t) \, d t = 0$. It remains to sum over $1 \leq l \leq m_k$:
$$\sum_{l=1}^{m_k} \big|\ln|z-\xi_{k,l}| - \int_{a_{k,l}}^{b_{k,l}} \ln |z-t|  V_k(t) dt \big| \lesssim \sum_{l=1}^{m_k} \int_{a_{k,l}}^{b_{k,l}} {{dt} \over {|z-t|^2}} \leq \int_{\R} {{dt} \over {(x-t)^2+y^2}} \lesssim {{1} \over {y}},$$
where the second inequality follows from $b_{k,l} = a_{k,l+1}$. The same argument applied to $[j_{k,l} - 1/2, j_{k,l} + 1/2]$ with $V_k(t) \mapsto 1$ and $\xi_{k,l} \mapsto j_{k,l}$ gives the second inequality of \eqref{eq:Gint}. \par

{\bf Step 2}. Second we prove
$$\bigg| \ln |G_k(z)| - \ln \left|{{\sin \pi z } \over {\pi}}\right| - \int_{\mathbb R} \ln |z - t| (V_k(t) - 1) \, d t \bigg| \lesssim {{1} \over {y}}.$$

In the notation above we have
$$\bigcup_{l = 1}^{m_k} [a_{k,l}, b_{k,l}] = \bigcup_{l = 1}^{m_k} [j_{k,l} - 1/2, j_{k,l} + 1/2] = [-4 \cdot 3^{k + 1} - 1/2, 4 \cdot 3^{k + 1} + 1/2].$$
Applying \eqref{eq:Gint} we obtain
\begin{equation}\label{eq:poissonfin}
\bigg| \sum_{l = 1}^{m_k} \ln \left|\frac{z - \xi_{k,l}}{z - j_{k,l}}\right| - \int_{\R} \ln |z - t| (V_k(t) - 1) \, d t \bigg| \lesssim  {{1} \over {y}},
\end{equation}
where the integral extends to the whole real line because $V_k - 1$ vanishes outside the common interval. The definition of $G_k$ therefore gives the claim.\par

{\bf Step 3}. Third we prove $\displaystyle \int_{\mathbb R} \ln |z - t| g_k(t) \, d t = P_{y + 4} \varphi_k(x)$. The Poisson formula for the logarithm states that
$$\frac{1}{\pi} \int_{\mathbb R} \frac{c \ln |x + i y - t|}{(t - s)^2 + c^2} \, d t = \ln |x + i (y + c) - s| \quad \text{for every } c > 0.$$
Note that for $c$ in a neighbourhood of $4$, the differentiated kernel decays quadratically in $t$. Hence, one can differentiate this identity with respect to $c$, multiply by $\varphi_k(s)/\pi$, integrate with respect to $s$, and set $c = 4$. By the definition of $g_k$ this gives the claim.\par 

{\bf Step 4}. It remains to show that replacing $V_k - 1$ by $g_k$ in \eqref{eq:poissonfin} gives a bounded error. Recall that 
$$g_k(t) = \frac{1}{\pi^2} \int_{\mathbb R} \varphi_k(s) \frac{(t - s)^2 - 16}{((t - s)^2 + 16)^2} \, d s \quad \text{and} \quad \int_{\mathbb R} g_k(t) \, d t = 0.$$
Corollary~\ref{cor:phiint} implies
$$|g_k(t)| \lesssim \frac{1}{t^2} \int_{\mathbb R} |\varphi_k(s)| \, d s \lesssim \frac{3^{k + 1}}{t^2} \quad \text{for every } |t| \geq 4 \cdot 3^{k + 1}.$$
For $|t| \leq 4 \cdot 3^{k + 1}$ we therefore obtain:
$$|V_k(t) - 1 - g_k(t)| =  \left| - \frac{1}{8 \cdot 3^{k + 1}} \int_{|s| < 4 \cdot 3^{k + 1}} g_k(s) \, d s \right| \leq \frac{1}{8 \cdot 3^{k + 1}} \int_{|s| > 4 \cdot 3^{k + 1}} |g_k(s)| \, d s $$
$$\lesssim   \frac{3^{k+1}}{8 \cdot 3^{k + 1}} \int_{|s| > 4 \cdot 3^{k + 1}} {{ds} \over {s^2}} = {{1} \over {16 \cdot 3^{k+1}}} .$$
For $|t| > 4 \cdot 3^{k + 1}$, the same difference equals $-g_k(t)$. Combining these observations we obtain
$$|V_k(t) - 1 - g_k(t)| \lesssim \frac{3^{k + 1}}{t^2 + 3^{2(k + 1)}}.$$
Since $\int_{\R} (V_k(t)-1-g_k(t)) \ dt = 0$, we have
$$\left| \int_{\R} \ln |z-t| (V_k(t) -1 - g_k(t)) \ dt \right| = \left| \int_{\R} \ln {{|z-t|} \over {|z+i3^{k+1}|}} (V_k(t) - 1 - g_k(t)) \ dt \right| $$
$$\lesssim \int_{\R} \left| \ln {{|z-t|} \over {|z+i3^{k+1}|}} \right| {{3^{k+1}} \over {t^2 + 3^{2(k+1)}}} \ dt.$$
It remains to bound the last integral uniformly in $k, x$ and $y$. The logarithmic Poisson formula with $s = 0$ and $c = 3^{k + 1}$ gives
$$\int_{\mathbb R} \ln \frac{|z - t|}{|z + i 3^{k + 1}|} \frac{3^{k + 1}}{t^2 + 3^{2(k + 1)}} \, d t = 0.$$
Note that the integral is absolutely convergent, since $y \ge 1$. Furthermore, $|z + i 3^{k + 1}| \geq \max\{|z|, 3^{k + 1}\}$, so the triangle inequality gives
$$\ln \frac{|z - t|}{|z + i 3^{k + 1}|} \leq \ln \left(1 + \frac{|t|}{3^{k + 1}}\right).$$
Now we obtain
$$\int_{\mathbb R} \left|\ln \frac{|z - t|}{|z + i 3^{k + 1}|}\right| \frac{3^{k + 1}}{t^2 + 3^{2(k + 1)}} \, d t = 2 \int_{\R} \max\left\{\ln {{|z-t|} \over {|z+i3^{k+1}|}}, 0 \right\} {{3^{k+1}} \over {t^2 + 3^{2(k+1)}}} \ dt $$
$$\leq 2 \int_{\R} \ln \left( 1+ {{|t|} \over {3^{k+1}}} \right) {{3^{k+1}} \over {t^2 + 3^{2(k+1)}}} \ dt  = 2 \int_{\mathbb R} \frac{\ln(1 + |u|)}{1 + u^2} \, d u < \infty,$$
where the last equality follows from the substitution $t \mapsto 3^{k+1} u$. Consequently,
$$\left|\int_{\mathbb R} \ln |z - t| (V_k(t) - 1 - g_k(t)) \, d t\right| \lesssim 1.$$

Combining this estimate with the identities above and using uniform boundedness of the function $\ln |\sin \pi (x + i y) / \pi| - \pi y $ for $y \geq 1$, we obtain the desired estimate. 

\end{proof}

\subsection{The final system is a Schauder basis}\label{sect:schauder}

The coordinate projections onto the first $l$ inserted exponentials are $S_k D_l M_k^{-1} F_k$ for $0 \leq l \leq m_k$, by Lemma~\ref{lem:algprefixes}. It is sufficient to show that their norms are bounded independently of $k$ and $l$. Once this bound is established, choose an absolute constant $C_1 \geq 1$ which bounds all these projection norms. Lemma~\ref{lem:algprefixes} gives
\begin{equation}\label{eq:Schauder}
D_l \Gram_k D_l \preceq 2 \pi C_1^2 M_k^* M_k
\end{equation}
and in particular $\|S_k M_k^{-1}\| \leq \sqrt{2 \pi} C_1$ by taking $l = m_k$. Together with the invertibility of $M_k$ established in Lemma~\ref{lem:density}, these bounds verify all the conditions needed in the proof of Theorem~\ref{th:algorithm} for the Schauder basis property. \par

The desired uniform bound for the projection norms is equivalent to
\begin{equation}\label{eq:Sch}
\left|\sum_{j = 1}^{l} (M_k^{-1} F_k f)_j \langle e_{\xi_{k,j}}, g\rangle\right| \lesssim \|f\|_2 \|g\|_2 \quad \text{for every } f, g \in L^2(-\pi, \pi).
\end{equation}
The case $l = 0$ is obvious. Now suppose $l >0$ and note that every finite Fourier projection $Q_E$ defined in \eqref{def:Q} has norm at most $1$. It is therefore enough to compare the sum above with $\langle Q_E f, g\rangle$ for a suitable finite set $E \subset \mathbb Z$, and then bound their difference by $\const \cdot \|f\|_2 \|g\|_2$. \par 
We will express this difference as a contour integral using the generating function $G_k$. We now introduce some notation. For $x \in \mathbb R$ and $y \geq 0$ put
\begin{equation}\label{eq:W}
W_k(x, y) = \min\left\{1, \left(\frac{1 + \dist(x, A_k) + y}{3^{k + 1}}\right)^{\frac{1}{2}}\right\}.
\end{equation}
For $f, g \in L^2(-\pi, \pi)$ define
$$b_k(z; f, g) = -\frac{\langle e_z, g\rangle}{G_k(z)} \sum_{n = 1}^{m_k} \frac{G_k(j_{k,n}) (F_k f)_n}{j_{k,n} - z}$$
with the convention that if $G_k(j_{k,n}) = 0$, then the corresponding summand equals zero. The following two technical lemmas provide the estimates and the contour identity needed for the argument.

\begin{lemma}\label{lem:int_estimates}
For every $x \in \mathbb R$, $f, g \in L^2(-\pi, \pi)$ and $y \geq 1$ we have
$$|b_k(x \pm i y; f, g)| \lesssim \frac{\|f\|_2 \|g\|_2}{y} \quad \text{and} \quad \int_1^\infty |b_k(x \pm i y; f, g)| \, d y \lesssim \|f\|_2 \|g\|_2.$$
Moreover, if $\dist(x, \mathcal Z_k) \geq 1/4$ and $\dist(x, \mathbb Z) \geq 1/8$, then $|b_k(x + i y; f, g)| \lesssim \|f\|_2 \|g\|_2$ for every $|y| \le 1$.
\end{lemma}

\begin{lemma}\label{lem:residues}
Let $1 \leq l \leq m_k$ and fix $a < b$ such that $a, b \notin \mathcal Z_k \cup \Z$. Assume that $\mathcal Z_k \cap (a, b) = \{\xi_{k,1}, \ldots, \xi_{k,l}\}$.
For $Y > 0$ let $\mathcal R = [a, b] + i [-Y, Y]$ with positively oriented boundary. Then
$$\sum_{j = 1}^{l} (M_k^{-1} F_k f)_j \langle e_{\xi_{k,j}}, g\rangle = \sum_{\substack{1 \leq n \leq m_k \\ a < j_{k,n} < b}} (F_k f)_n \langle e_{j_{k,n}}, g\rangle - \frac{1}{2 \pi i} \int_{\partial \mathcal R} b_k(z; f, g) \, d z.$$
\end{lemma}

We first deduce estimate \eqref{eq:Sch} from these lemmas. Fix $1 \leq l \leq m_k$. The inserted frequencies form a consecutive block in $\mathcal Z_k$. Start with the midpoints of the gaps immediately before $\xi_{k,1}$ and immediately after $\xi_{k,l}$. Each midpoint has distance at least $3/7$ from $\mathcal Z_k$ (see Lemma~\ref{lem:density}).\par

If necessary, move each midpoint by at most $1/8$ so that its distance from $\Z$ becomes at least $1/8$. Denote the resulting points by $a$ and $b$. They remain in the same gaps and have distance at least $1/4$ from $\mathcal Z_k$. Thus they satisfy the hypotheses of both lemmas.\par

Take $ Y = \max\{2, b - a\}$ and $\mathcal R = [a, b] + i [-Y, Y]$ with positively oriented boundary. We need to bound the integral from Lemma~\ref{lem:residues}. On each horizontal side, parameterized by $z = x \pm i Y$ for $a \leq x \leq b$, the pointwise estimate in Lemma~\ref{lem:int_estimates} gives
$$\int_a^b |b_k(x \pm i Y; f, g)| \, d x \lesssim \frac{b - a}{Y} \|f\|_2 \|g\|_2 \leq \|f\|_2 \|g\|_2.$$

Set $s \in \{a,b\}$ and split the corresponding vertical side into three segments: the upper and lower segments with $z = s \pm i y$, where $1 \leq y \leq Y$, and the middle segment with $z = s + i y$, where $-1 \leq y \leq 1$. The integral estimate in Lemma~\ref{lem:int_estimates} gives
$$\int_1^Y |b_k(s \pm i y; f, g)| \, d y \leq \int_1^\infty |b_k(s \pm i y; f, g)| \, d y \lesssim \|f\|_2 \|g\|_2.$$
The distance conditions on $s$ allow us to apply the last assertion of Lemma~\ref{lem:int_estimates} to the middle segment:
$$\int_{-1}^1 |b_k(s + i y; f, g)| \, d y \lesssim \|f\|_2 \|g\|_2.$$
Adding the estimates for all boundary segments, we obtain
$$\left|\int_{\partial \mathcal R} b_k(z; f, g) \, d z\right| \leq \int_{\partial \mathcal R} |b_k(z; f, g)| \, |d z| \lesssim \|f\|_2 \|g\|_2.$$

The finite sum on the right-hand side of the identity in Lemma~\ref{lem:residues} equals $\langle Q_E f, g\rangle$ for $E = \{j_{k,n} \mid a < j_{k,n} < b \text{ and } 1 \leq n \leq m_k \}$, since $(F_k f)_n = \widehat f(j_{k,n})$. Since $\|Q_E\|\leq 1$, we obtain $| \langle Q_E f, g\rangle|\lesssim \|f\|_2 \|g\|_2$. Using Lemma~\ref{lem:residues}, we therefore have
$$\left| \left\langle S_k D_l M_k^{-1} F_k f, g \right\rangle \right| \leq \left| \langle Q_E f, g\rangle \right| + \frac{1}{2 \pi} \left| \int_{\partial \mathcal R} b_k(z; f, g) \, d z \right| \lesssim \|f\|_2 \|g\|_2.$$
This proves \eqref{eq:Sch}. It remains to prove the two technical lemmas.\par

\begin{proof}[Proof of Lemma~\ref{lem:int_estimates}]
Lemma~\ref{lem:genprofile} and the identity $G_k(\overline z) = \overline{G_k(z)}$ give
\begin{equation}\label{eq:Gy}
|G_k(x \pm i y)|^2 \asymp e^{2 \pi y} W_k(x, y) \quad \text{for every } y \geq 1,
\end{equation}
where $W_k$ is defined in \eqref{eq:W}. We also need the corresponding upper bound on the real axis. The product formula for $\sin \pi z$ gives
\begin{equation} \label{eq:Gprod}
\frac{|G_k(x + i y)|^2}{|G_k(x)|^2} = \prod_{t \in \mathcal Z_k} \left(1 + \frac{y^2}{(x - t)^2}\right) \quad \text{for every } x \in \R \setminus \mathcal Z_k \text{ and } y \in \R,
\end{equation}
since the identity extends by continuity to deleted integer points which are not zeros of $G_k$. Every factor in \eqref{eq:Gprod} is at least $1$. Taking $y = 1$ and using \eqref{eq:Gy} with the estimate $W_k(x, 1) \leq \sqrt{2} W_k(x, 0)$, we obtain 
\begin{equation}\label{eq:Greal}
|G_k(x)|^2 \lesssim W_k(x, 0) \quad \text{for every } x \in \R.
\end{equation} 

Fix $x \in \R$. We first prove the lemma for $b_k(x+iy; f, g)$ for $y \geq 1$. The definition of $W_k$ gives
\begin{equation}\label{eq:Wineq}
W_k(t, 0) \leq W_k(x, |t - x|) \quad \text{and} \quad W_k(x, s) \leq \max\left\{1, \sqrt{{{s} \over {y}}} \right\} W_k(x, y) \quad \text{for every } s, y \geq 1.
\end{equation}
Indeed, the first inequality follows from $\dist(t, A_k) \leq \dist(x, A_k) + |t - x|$. The second one follows from monotonicity when $s \leq y$ and from the estimate
$$1 + \dist(x, A_k) + s \leq {{s} \over{y}} \left(1 + \dist(x, A_k) + y\right)$$
when $s \geq y$.  Put $p_{k,n} = 1 + |j_{k,n} - x|$ for $1 \leq n \leq m_k$. For every $y \geq 1$ put
$$v_{k,n}(y) = \frac{G_k(j_{k,n})}{\sqrt{W_k(x, y)} (j_{k,n} - x - i y)}.$$
Combining \eqref{eq:Greal} with $|j_{k,n} - x - i y| \asymp p_{k,n} + y$, we obtain
\begin{equation}\label{eq:vkn}
|v_{k,n}(y)| \lesssim \frac{\sqrt{W_k(j_{k,n}, 0)}}{\sqrt{W_k(x, y)} (p_{k,n} + y)} \leq \frac{\sqrt{W_k(x, |j_{k,n} - x|)}}{\sqrt{W_k(x, y)} (p_{k,n} + y)} \leq {{\sqrt{W_k(x,p_{k,n})}} \over {\sqrt{W_k(x,y)} (p_{k,n}+y)}} \leq \frac{\max\left\{1, (p_{k,n}/y)^{1/4}\right\}}{p_{k,n} + y},
\end{equation}
where the third inequality follows from monotonicity of $W_k$ and the last one follows from \eqref{eq:Wineq}. Thus $|v_{k,n}(y)| \lesssim 1/y$ when $p_{k,n} \leq y$, and $|v_{k,n}(y)| \lesssim y^{-1/4} p_{k,n}^{-3/4}$ when $p_{k,n} > y$.

Observe that every interval of length $1$ contains at most four of the numbers $p_{k,n}$, counted with multiplicity, since the $j_{k,n}$ are distinct integers. Using \eqref{eq:vkn}, we then obtain
\begin{equation}\label{eq:vsum}
\sum_{n = 1}^{m_k} \frac{|v_{k, n}(y)|}{\sqrt{p_{k,n}}} \lesssim \frac{1}{y} \sum_{p_{k,n} \leq y} {{1} \over {\sqrt{p_{k, n}}}} + {{1} \over {y^{1/4}}} \sum_{p_{k,n} > y} p_{k, n}^{-5/4} \lesssim \frac{1}{\sqrt y}.
\end{equation} 
The same bound from \eqref{eq:vkn} also gives
\begin{equation}\label{eq:vmax}
\max_n \sqrt{p_{k,n}} |v_{k,n}(y)| \lesssim {{1} \over {\sqrt{y}}}.
\end{equation} 
Applying \eqref{eq:vkn} separately on $1 \leq y \leq p_{k,n}$ and $y \geq p_{k,n}$, we obtain
\begin{equation}\label{eq:vintegral}
\int_1^\infty \frac{|v_{k,n}(y)|}{\sqrt y} \, d y \lesssim p_{k,n}^{-3/4} \int_1^{p_{k,n}} y^{-3/4} \, d y + \int_{p_{k,n}}^\infty y^{-3/2} \, d y \lesssim {{1} \over {\sqrt{p_{k,n}}}}.
\end{equation} 

We first use these coefficient estimates to bound the finite sum in $b_k$ pointwise in $y$. We have
\begin{equation}\label{eq:vsquare}
\sum_{n = 1}^{m_k} |v_{k,n}(y)|^2 \leq \left(\max_{1 \leq n \leq m_k} \sqrt{p_{k,n}} |v_{k,n}(y)|\right) \sum_{n = 1}^{m_k} \frac{|v_{k,n}(y)|}{\sqrt{p_{k,n}}} \lesssim \frac{1}{y}
\end{equation}
by \eqref{eq:vsum} and \eqref{eq:vmax}. Using both Cauchy--Schwarz and Bessel inequalities, we obtain
\begin{equation}\label{eq:sum-pointwise}
\left|\sum_{n = 1}^{m_k} v_{k,n}(y) (F_k f)_n\right|^2 \leq \left(\sum_{n = 1}^{m_k} |v_{k,n}(y)|^2\right) \sum_{n = 1}^{m_k} |(F_k f)_n|^2 \lesssim \frac{1}{y} \cdot \|f\|_2^2.
\end{equation}

We next bound the $L^2(1,\infty)$-norm of the same sum. One can obtain
$$\left|\sum_{n = 1}^{m_k} v_{k,n}(y) (F_k f)_n\right|^2 \leq \left(\sum_{n = 1}^{m_k} \frac{|v_{k,n}(y)|}{\sqrt{p_{k,n}}}\right) \left(\sum_{n = 1}^{m_k} \sqrt{p_{k,n}} |v_{k,n}(y)| |(F_k f)_n|^2\right).$$
Apply \eqref{eq:vsum} to the first factor and integrate in $y$. Estimate \eqref{eq:vintegral} then gives
\begin{equation}\label{eq:sum-integral}
\int_1^\infty \left|\sum_{n = 1}^{m_k} v_{k,n}(y) (F_k f)_n\right|^2 \, d y \lesssim \sum_{n = 1}^{m_k} \sqrt{p_{k,n}} |(F_k f)_n|^2 \int_1^\infty \frac{|v_{k,n}(y)|}{\sqrt y} \, d y \lesssim \|f\|_2^2.
\end{equation}

We now obtain the corresponding estimates for the exponential factor in $b_k$. One has
\begin{equation}\label{eq:exp-pointwise}
e^{-\pi y} |\langle e_{x + i y}, g\rangle| \leq \|g\|_2 \cdot e^{-\pi y} \sqrt{\int_{-\pi}^{\pi} e^{-2yt} \ dt} \lesssim \frac{\|g\|_2}{\sqrt y} \quad \text{for every } y \geq 1.
\end{equation}
Using $\int_0^\infty e^{-y t} t^{-1/2} \, d t \lesssim y^{-1/2}$ we have
$$e^{-2 \pi y} |\langle e_{x + i y}, g\rangle|^2 = \left| \int_0^{2 \pi} e^{-y t} e^{i x(t - \pi)} \overline{g(t - \pi)} \, d t  \right|^2 \lesssim \frac{1}{\sqrt y} \int_0^{2 \pi} e^{-y t} \sqrt t \, |g(t - \pi)|^2 \, d t.$$
Integrate in $y$ and use $\int_1^\infty e^{-y t} y^{-1/2} \, d y \lesssim t^{-1/2}$ for $t > 0$:
\begin{equation}\label{eq:exp-integral}
\int_1^\infty e^{-2 \pi y} |\langle e_{x + i y}, g\rangle|^2 \, d y \lesssim \int_0^{2 \pi} \sqrt t \, |g(t - \pi)|^2 \left(\int_1^\infty \frac{e^{-y t}}{\sqrt y} \, d y\right) d t \lesssim \|g\|_2^2.
\end{equation}
The definition of $v_{k,n}$ and \eqref{eq:Gy} give
$$|b_k(x + i y; f, g)| \lesssim e^{-\pi y} |\langle e_{x + i y}, g\rangle| \left|\sum_{n = 1}^{m_k} v_{k,n}(y) (F_k f)_n\right|.$$
Now we combine \eqref{eq:sum-pointwise} with \eqref{eq:exp-pointwise} to prove
$$|b_k(x + i y; f, g)| \lesssim {{ \|f\|_2 \|g\|_2} \over {y}}.$$
Applying Cauchy--Schwarz in $y$ and then using \eqref{eq:sum-integral} and \eqref{eq:exp-integral} proves 
$$\int_1^\infty |b_k(x + i y; f, g)| \, d y \lesssim \|f\|_2 \|g\|_2.$$

The estimates for $b_k(x - i y; f, g)$ with $y \geq 1$ follow by the same argument. Indeed, replacing $j_{k,n} - x - i y \mapsto j_{k,n} - x + i y$ leaves all coefficient bounds unchanged, and we have $\langle e_{x - i y}, g\rangle = \langle e_{-x + i y}, g(-\cdot)\rangle$.\par 

Assume that $\dist(x, \mathcal Z_k) \geq 1/4$ and $\dist(x, \mathbb Z) \geq 1/8$. Now $\sum_{t \in \mathcal Z_k} |x - t|^{-2} \lesssim 1$ because of the first assumption and the fact that $\mathcal{Z}_k$ is uniformly discrete. Taking logarithms in \eqref{eq:Gprod} at $y = 1$ we obtain $|G_k(x + i)|^2 \lesssim |G_k(x)|^2$. Each factor in \eqref{eq:Gprod} is nondecreasing in $|y|$. Together with \eqref{eq:Gy} at $y = 1$ this gives
\begin{equation}\label{eq:Gstrip}
|G_k(x + i y)|^2 \asymp |G_k(x + i)|^2 \asymp W_k(x, 1) \quad \text{for every } |y| \leq 1.
\end{equation}

By assumption $|j_{k,n} - x| \geq 1/8$, and hence $|j_{k,n} - x - i| \lesssim |j_{k,n} - x - i y|$ uniformly for $|y| \leq 1$. By the definition of $v_{k,n}(1)$ we obtain
$$\left|\sum_{n = 1}^{m_k} \frac{G_k(j_{k,n}) (F_k f)_n}{j_{k,n} - x - i y}\right| \lesssim \sqrt{W_k(x, 1)} \left(\sum_{n = 1}^{m_k} |v_{k,n}(1)|^2\right)^{1/2} \left(\sum_{n = 1}^{m_k} |(F_k f)_n|^2\right)^{1/2}.$$
Applying \eqref{eq:vsquare} at $y = 1$ and Bessel's inequality gives the bound $\const \cdot \sqrt{W_k(x, 1)} \|f\|_2$ for the right-hand side. Also, $|\langle e_{x + i y}, g\rangle| \lesssim \|g\|_2$ for $|y| \leq 1$, since the norms $\|e_{x + i y}\|_2$ are uniformly bounded on this strip. Substituting these estimates into the definition of $b_k$ and dividing by the lower bound in \eqref{eq:Gstrip}, we obtain $|b_k(x + i y; f, g)| \lesssim \|f\|_2 \|g\|_2$.

\end{proof}

\begin{proof}[Proof of Lemma~\ref{lem:residues}]
By the definition of $G_k$ we have
\begin{equation} \label{eq:gkg} \frac{G_k(z)}{G_k'(\xi_{k,j}) (z - \xi_{k,j})} = \frac{\sin(\pi z)}{\pi} L_{k,j}(z),
\end{equation}
where $L_{k,j}$ is some rational function with possible simple poles at $j_{k,1}, \ldots, j_{k,m_k}$. For fixed $k$ and $j$ it satisfies $|L_{k,j}(z)| \leq C(k,j) 1/|z|$ as $z \to \infty$. Hence, its partial fraction decomposition has no polynomial term:
$$L_{k,j}(z) = \sum_{n = 1}^{m_k} \frac{\operatorname{Res}_{w = j_{k,n}} L_{k,j}(w)}{z - j_{k,n}}.$$
Multiplying by $\sin(\pi z)/\pi$ expresses the entire function from \eqref{eq:gkg} as a linear combination of $\sinc(z - j_{k,n})$. Its values are $1$ at $\xi_{k,j}$ and $0$ at every other inserted frequency. By the definition of $M_k$, the coefficients form the $j$-th row of $M_k^{-1}$. Thus
\begin{equation}\label{eq:cardinal-expansion}
\frac{G_k(z)}{G_k'(\xi_{k,j}) (z - \xi_{k,j})} = \sum_{n = 1}^{m_k} (M_k^{-1})_{j n} \sinc(z - j_{k,n}).
\end{equation}
Each quotient ${{G_k(z) - G_k(j_{k,i})} \over {z - j_{k,i}}}$ extends to an entire function with value $G_k'(j_{k,i})$ at $j_{k,i}$. Evaluating \eqref{eq:cardinal-expansion} at $z = j_{k,n}$ gives 
$$(M_k^{-1})_{j n} = {{G_k(j_{k,n})} \over {G_k'(\xi_{k,j}) (j_{k,n} - \xi_{k,j})}},$$
when $j_{k,n} \ne \xi_{k,j}$. Since $G_k(\xi_{k,j}) = 0$, this can be rewritten as
\begin{equation}\label{eq:rescoef}
\left.\frac{G_k(z) - G_k(j_{k,n})}{z - j_{k,n}}\right|_{z = \xi_{k,j}} = G_k'(\xi_{k,j}) (M_k^{-1})_{j n}.
\end{equation}
This also holds when $j_{k,n} = \xi_{k,j}$: here $(M_k^{-1})_{j n} = 1$ and both sides of \eqref{eq:rescoef} equal $G_k'(\xi_{k,j})$.

We now rewrite the definition of $b_k$ as
\begin{equation}\label{eq:resdec}
b_k(z; f, g) = \langle e_z, g\rangle \sum_{n = 1}^{m_k} \frac{(F_k f)_n}{z - j_{k,n}}  - \frac{\langle e_z, g\rangle}{G_k(z)} \sum_{n = 1}^{m_k} (F_k f)_n \frac{G_k(z) - G_k(j_{k,n})}{z - j_{k,n}}.
\end{equation}
The function $z \mapsto \langle e_z, g\rangle$ is entire. Thus the first summand in \eqref{eq:resdec} has residue $(F_k f)_n \langle e_{j_{k,n}}, g\rangle$ at each deleted integer $j_{k,n}$. The residue of the second summand at $\xi_{k,j}$ is $-(M_k^{-1} F_k f)_j \langle e_{\xi_{k,j}}, g\rangle$ by \eqref{eq:rescoef}.

By the assumptions of the lemma the only zeros of $G_k$ inside $\mathcal R$ are $\xi_{k,1}, \ldots, \xi_{k,l}$. Neither term in \eqref{eq:resdec} has a pole on the boundary. Applying the residue theorem to each term gives the identity in the lemma. If a deleted integer coincides with an inserted frequency, its two residue contributions are automatically combined.
\end{proof}

\subsection{The final system is not a Riesz basis}

It remains to verify the eigenvalue condition in Theorem~\ref{th:exist}. We have $\|S_k c\|_2^2 \geq \lambda_{\min}(\Gram_k) \|c\|_2^2$ for every $c \in \CC^{m_k}$. Applying this to $M_k^{-1} c$ and taking the supremum over $\|c\|_2=1$, we obtain
$$\lambda_{\min}(\Gram_k) \|M_k^{-1}\|^2 \leq \|S_k M_k^{-1}\|^2 \leq 2 \pi C_1^2$$
by \eqref{eq:Schauder}. It remains to obtain a lower bound for the norm of $M_k^{-1}$.

\begin{lemma} We have $\|M_k^{-1}\|^2 \gtrsim 2^{k + 1} 3^{-\frac{k + 1}{2}}$.
\end{lemma}

\begin{proof}
For each point $t \in A_k^+$ choose a nearest zero $\xi_{k,l}$ of $G_k$. Define $E_k$ to be the set of indices of the selected zeros. We have $|t - \xi_{k,l}| = \dist(t, \mathcal Z_k) \leq 3/5$ by Lemma~\ref{lem:density}. Since $A_k^+$ has $2^{k + 1}$ points in $[0, 3^{k + 1} - 1]$ separated by at least $2$, the selected zeros are distinct and they belong to the inserted block. Obviously they satisfy $|\xi_{k,l}| < 3^{k + 1}$.\par

Using \eqref{eq:Gprod} with $y = 1$ we obtain
$$|G_k(x+i)|^2 \geq |G_k(x)|^2 \left( 1+ {{1} \over {(x-\xi_{k,l})^2}}\right) \geq {{|G_k(x)|^2} \over {(x-\xi_{k,l})^2}}.$$
Letting $x \to \xi_{k,l}$ we get $|G_k (\xi_{k,l}+i)|^2 \geq |G_k' (\xi_{k,l})|^2$. But then \eqref{eq:W} and \eqref{eq:Gy} give
$$|G_k(\xi_{k,l}+i)|^2 \lesssim W_k (\xi_{k,l},1) \leq \left( {{2+3/5} \over {3^{k+1}}} \right)^{1/2} \lesssim {{1} \over {3^{\frac{k+1}{2}}}}.$$
Combining these bounds, we have
\begin{equation}\label{eq:selected-derivatives}
|G_k'(\xi_{k,l})| \leq |G_k(\xi_{k,l} + i)| \lesssim 3^{-{{k + 1} \over {4}}} \quad \text{for } l \in E_k.
\end{equation}
Put $z_k = i 3^{k + 1}$. Note that $W_k(0, 3^{k + 1}) = 1$. Now Bessel's inequality and \eqref{eq:Gy} give
$$\sum_{n = 1}^{m_k} |\sinc(z_k - j_{k,n})|^2 = {{1} \over {4\pi^2}} \sum_{n = 1}^{m_k} | \langle e_{z_k}, e_{j_{k,n}}\rangle|^2 \leq \frac{\|e_{z_k}\|_2^2}{2 \pi} \lesssim {{e^{2\pi \cdot 3^{k+1}} W_k(0,3^{k+1})} \over {3^{k+1}}} \asymp \frac{|G_k(z_k)|^2}{3^{k + 1}}.$$
Combining this with \eqref{eq:cardinal-expansion} yields
$$|G_k(z_k)|^2 \sum_{l \in E_k} \frac{1}{|G_k'(\xi_{k,l})|^2 |z_k - \xi_{k,l}|^2} = \sum_{l \in E_k} \left|\sum_{n = 1}^{m_k} (M_k^{-1})_{l n} \sinc(z_k - j_{k,n})\right|^2 \leq \sum_{l =1}^{m_k} \left|\sum_{n = 1}^{m_k} (M_k^{-1})_{l n} \sinc(z_k - j_{k,n})\right|^2 $$
$$\leq \|M_k^{-1}\|^2 \sum_{n = 1}^{m_k} |\sinc(z_k - j_{k,n})|^2 \lesssim {{\|M_k^{-1}\|^2 |G_k(z_k)|^2} \over {3^{k+1}}}.$$
Since all zeros of $G_k$ are real, we have $G_k(z_k) \ne 0$. So we may cancel $|G_k(z_k)|^2$ and use the estimate \eqref{eq:selected-derivatives}:
$$\sum_{l \in E_k} \frac{1}{|G_k'(\xi_{k,l})|^2 |z_k - \xi_{k,l}|^2} \gtrsim \sum_{l \in E_k} \frac{3^{(k + 1)/2}}{3^{2(k + 1)}} = {{1} \over {3^{\frac{3(k + 1)}{2}}}} \cdot \# E_k = {{2^{k+1}} \over {3^{\frac{3(k + 1)}{2}}}},$$
since $|z_k - \xi_{k,l}|^2 = |i3^{k+1} - \xi_{k,l}|^2 \leq 2 \cdot 3^{2(k + 1)}$.  Finally, we obtain
$$\|M_k^{-1}\|^2 \gtrsim 3^{k + 1} \sum_{l \in E_k} \frac{1}{|G_k'(\xi_{k,l})|^2 |z_k - \xi_{k,l}|^2} \gtrsim 2^{k + 1} 3^{-\frac{k + 1}{2}}.$$
\end{proof}

Hence, there exists an absolute constant $C_2 \geq 1$ such that for every $k \in \N$ we have
\begin{equation}\label{eq:c2}
\lambda_{\min}(\Gram_k) \leq C_2^2 \left(\frac{3}{4}\right)^{{{k + 1} \over {2}}}.
\end{equation}
This completes the proof of Theorem \ref{th:exist}.

\section{Acknowledgments and AI use statement}
The work is supported by RSF grant 26-11-00157. The author wishes to express his deep gratitude to Anton Tselischev, who introduced him to this problem in early 2026, and to Yurii Belov for helpful discussions of the underlying ideas. \par
The central idea of the manuscript, the Algorithm of finite replacements and the proof of its correctness were developed without the use of LLM. The ChatGPT Astra 6.0 was used to search for an initial example of such a basis satisfying ideas and constraints imposed by the algorithm. The basis presented in the article was obtained by the author through a substantial human-made rework of the initial example. The proofs of the technical lemmas presented in Section~\ref{sect:Tech} were considerably simplified compared with their original versions with the assistance of ChatGPT. The TeX-code for Figure 1 was generated by LLM.

\end{document}